%% file: main.tex
\documentclass[preprint,12pt]{elsarticle}

\usepackage{amsmath,amssymb,amsthm}
\usepackage{graphicx}
\usepackage{booktabs}

\usepackage{natbib}

\usepackage{enumitem}
\usepackage{subcaption}
\usepackage{tikz}
\usetikzlibrary{calc,backgrounds,arrows.meta,decorations.pathreplacing}
\usepackage{pgfplots}
\pgfplotsset{compat=1.18}
\usepackage{pgfplotstable}
\usepackage{siunitx}
\usepackage{threeparttable}
\usepackage{csvsimple}
\usepackage[section]{placeins}

\input{macros}

\begin{document}

\journal{Computers \& Industrial Engineering}

\begin{frontmatter}

\title{Exact formulations for rectangular-warehouse single-picker routing with scattered storage in single-block and two-block layouts}

\author[inst1]{George Dunn\corref{cor1}}

\cortext[cor1]{Corresponding author.}
\ead{george.dunn@uon.edu.au}

\address[inst1]{School of Computer and Information Sciences, The University of Newcastle, University Drive, Callaghan NSW 2308, Australia}

\begin{abstract}
Order picking travel dominates much of warehouse effort, and exact routing is especially valuable when storage is scattered so pick locations are not fixed in advance.
We address the single picker routing problem (SPRP) and its scattered-storage variant (SPRP-SS) in single-block and two-block rectangular warehouses.
We propose two mixed-integer linear programming formulations that exploit structural properties of optimal tours to simplify connectivity modelling and remove redundant edge configurations: a \emph{Configuration Connectivity} model tailored to single-block layouts and an \emph{Edge Connectivity} model that extends to two-block layouts.
In extensive computational experiments on large randomly generated benchmark sets for single-block and two-block rectangular layouts, we compare these formulations against established MILP and network-flow baselines for SPRP and SPRP-SS and report computational gains tied to the structural restrictions.
The results support using compact, solver-based exact routing models in industrial settings where dynamic programming is cumbersome to integrate, particularly for SPRP-SS and for routing subproblems embedded in larger planning or warehouse-design optimizations.
\end{abstract}

\begin{keyword}
Order picking \sep Picker routing \sep Scattered storage \sep Mixed-integer linear programming \sep Combinatorial optimization \sep Warehouse logistics
\end{keyword}

\end{frontmatter}

\section{Introduction}
\input{sections_introduction}

\section{Problem setting and notation}
\input{sections/section_problem_setting}

\section{Related work}
\label{sec:related}
\input{sections_related}

\section{Configuration Connectivity (CC) formulation (single-block)}
\input{sections/section_cc}

\section{Edge Connectivity (EC) formulation (single- and two-block)}
\input{sections/section_ec}

\section{Computational experiments}
\input{sections/section_results}

\section{Conclusions and implications}
\input{sections_conclusion}

\section*{Funding}
This work was supported by the Australian Government Research Training Program (RTP) Scholarship.

\section*{Data and Code Availability}
The numerical results were produced with author-written Python code and pseudorandom instance generators.
The code and scripts are maintained in a private repository for now and are available from the corresponding author on reasonable request.

\bibliographystyle{elsarticle-harv}
\bibliography{biblio}

\appendix
\input{sections/appendix_gs}

\end{document}

%% file: macros.tex
\newcommand{\SPRP}{\textsc{SPRP}}
\newcommand{\SPRPsS}{\textsc{SPRP-SS}}

\newcommand{\Jcal}{\mathcal{J}}
\newcommand{\Ical}{\mathcal{I}}
\newcommand{\Hcal}{\mathcal{H}}
\newcommand{\Ncal}{\mathbb{N}_0}
\newcommand{\xx}[1]{\underline{\underline{#1}}}   
\newcommand{\xb}[1]{\overline{\overline{#1}}}    
\newcommand{\xbx}[1]{\overline{\underline{#1}}}  
\newcommand{\xxb}[1]{\overline{\overline{\underline{\underline{#1}}}}}  
\newcommand{\xpass}[1]{x^{\vert}_{#1}}
\newcommand{\xppass}[1]{x^{\vert\vert}_{#1}}
\newcommand{\pxji}[2]{\underset{\scriptscriptstyle\wedge}{#1}_{#2}}
\newcommand{\qxji}[2]{\overset{\scriptscriptstyle\vee}{#1}_{#2}}
\newcommand{\pitop}[1]{\overset{\circ}{\pi}_{#1}}
\newcommand{\pibot}[1]{\underset{\circ}{\pi}_{#1}}
\newcommand{\cond}[2]{\underset{#1}{[}#2]}
\newcommand{\thetad}{\overset{\circ}{\theta}}

\newcommand{\Kcal}{\mathcal{K}}
\newcommand{\acbarx}[2]{\bar{x}_{#1#2}}
\newcommand{\acbarbarx}[2]{\bar{\bar{x}}_{#1#2}}
\newcommand{\acpass}[2]{x^{\vert}_{#1#2}}
\newcommand{\acr}[3]{\tau_{#1,#2,#3}}
\newcommand{\acp}[3]{\breve{p}_{#1(#2,#3)}}
\newcommand{\acz}[3]{\hat{z}_{#1(#2,#3)}}
\newcommand{\aceta}[1]{\eta_{#1}}
\newcommand{\acpi}[2]{\pi_{#1#2}}

\newcounter{modelctr}
\renewcommand{\themodelctr}{\arabic{modelctr}}

\newcommand{\EC}{%
  \setcounter{modelctr}{0}%
  \def\modelprefix{EC-}%
}

\newcommand{\CC}{%
  \setcounter{modelctr}{0}%
  \def\modelprefix{CC-}%
}

\newcommand{\modeltag}{%
  \refstepcounter{modelctr}%
  \tag{\modelprefix\themodelctr}%
}

%% file: sections_introduction.tex
Order picking (retrieving items for a customer or replenishment wave and returning to a depot) is a dominant driver of labour cost and throughput in many warehouses \citep{petersen1997evaluation,de2007design,tompkins2010facilities}.
A core tactical subproblem is the \emph{single picker routing problem} (SPRP): given a rectangular warehouse with parallel aisles, determine a minimum-distance tour that starts and ends at the depot and visits all required pick locations \citep{de2007design}.
The problem is a structured, $\mathcal{NP}$-hard variant of the traveling salesperson problem; while traversal-style heuristics remain common in practice \citep{hall1993distance}, exact optimisation is increasingly relevant when routing is embedded inside larger planning tasks (e.g., storage-policy evaluation, slotting, layout what-if analysis, and wave planning) that require solving many routing instances reliably.

Under \emph{scattered storage}, each SKU may be stored in multiple candidate locations, so the decision of \emph{where} to pick and \emph{how} to route are coupled.
This yields the SPRP with scattered storage (\SPRPsS{}), for which solver-based models are particularly attractive because they support integrated constraints and can be called as submodels within broader mixed-integer programs \citep{gu2007research,rouwenhorst2000warehouse,boysen202550}.
For these applications, the practical bottleneck is often not modelling expressiveness but \emph{formulation efficiency}: the number of variables and constraints, the strength of the linear relaxation, and avoidable symmetry can determine whether an exact routing model is usable at scale.

This paper develops \emph{structure-aware} MILP formulations for SPRP and \SPRPsS{} in \emph{single-block} and \emph{two-block} rectangular warehouses.
Recent exact formulations for these layouts sit between two useful benchmarks.
The edge-configuration MILP of \citet{goeke2021modeling} (GS) provides a compact, solver-friendly model and extends naturally to scattered storage, but it still represents configuration patterns that are unnecessary in rectangular layouts and therefore enlarge the search space.
At the other end, the network-flow formulation of \citet{hessler2024exact} (NF) models a shortest path on a dynamic-programming-derived state graph and has excellent performance in many regimes, but its state-graph construction and scaling behaviour differ markedly from configuration-based MILPs, which matters when routing is used as an embedded component rather than as a stand-alone algorithmic solver.

Our contribution is \emph{modelling}: we show how to encode rectangular-warehouse routing using smaller and stronger MILPs by aligning the decision variables and connectivity logic with regularities observed in optimal tours.
First, we derive a \emph{Configuration Connectivity} (CC) formulation for single-block layouts by trimming redundant configuration options from GS and by expressing connectivity using a reduced set of configuration decisions, which decreases model size and removes symmetry without changing the optimal objective value.
Second, we propose an \emph{Edge Connectivity} (EC) formulation that enforces connectivity directly on the warehouse graph through a unified core and layout-specific connectivity constraints, yielding a formulation family that covers both single-block and two-block layouts and supports \SPRPsS{} in the same framework.

\paragraph{Structural regularities exploited (informal)}
Optimal tours in rectangular layouts exhibit regular patterns: within an aisle, once the tour enters a subaisle it either traverses it once or performs a contiguous ``branch-and-pick'' to a furthest visited position before returning, and horizontal moves between adjacent aisles occur through cross-aisles in a way that can be summarised by a small set of configuration indicators.
These regularities motivate our variable choices and allow us to eliminate configuration patterns that do not improve feasibility or cost in the layouts studied.
Importantly, the MILP formulations themselves remain complete and verifiable from the constraints stated in this paper; prior structural analyses are used only as background to explain why the resulting formulations are compact and effective.
We cite the published elimination result in \citet{dunn2025double} and related manuscripts under review for additional intuition and refinements, but the formulations below do not require the reader to consult those sources to understand the model or assess correctness.

We benchmark CC and EC against GS and NF on large random instance sets for SPRP and \SPRPsS{} in single- and two-block layouts.
The results demonstrate that incorporating layout-specific structure at the formulation level can materially reduce model size and improve solution times, especially for \SPRPsS{} and for regimes where a compact MILP is preferable as a reusable building block in industrial decision-support models.
The remainder of the paper is organised as follows: Section~\ref{sec:related} reviews related exact approaches and formulations; Sections~4--5 present CC and EC; Section~6 reports computational results and regime-level guidance; and Section~7 concludes with implications for practice and extensions.

%% file: sections/section_problem_setting.tex
This section sets the context for the formulations that follow.
Up until now, we have only considered the standard SPRP,
where each requested \emph{stock keeping unit} (SKU)
is stored in a unique pick location.
Whilst this problem is still considered,
this paper will also address the more complex
\emph{scattered-storage} variation (SPRP-SS),
in which each SKU may be stored at multiple positions
so that the picker must choose both where to pick and which route to take.
The paper focuses on single-block and two-block rectangular layouts
but briefly relates these to the broader literature on multi-block
and mixed-shelves warehouses.


\subsection{Rectangular warehouse layouts}
\label{sec:layout}

We consider conventional rectangular warehouses comprised of parallel picking aisles
connected by horizontal cross-aisles.
In the simplest \emph{single-block} layout, the aisles run between two main cross-aisles.
Larger warehouses often include additional cross-aisles that partition the storage area into multiple blocks;
in this paper, we focus on the practically important \emph{two-block} case with exactly one additional cross-aisle.
Figure~\ref{fig:multiblock_warehouse} illustrates the block structure and relevant elements (aisles, cross-aisles, subaisles, and pick locations).

\input{figures/figure_multiblock_warehouse}

A layout is said to be \emph{rectangular} if all aisles are parallel,
all cross-aisles are parallel, aisles are perpendicular to cross-aisles,
and the spacing between consecutive aisles and between consecutive cross-aisles is constant.
Items are stored in \emph{subaisles} (the portion of an aisle between two cross-aisles), which are assumed narrow enough
that horizontal traversal within a subaisle is negligible.

\subsection{Graph representation}
\label{sec:graph}

To model routing decisions, the warehouse is represented as an undirected weighted graph
\(G = (V \cup P, E)\),
where \(V\) denotes aisle-intersection vertices and \(P\) denotes pick locations, including the depot.
Figure~\ref{fig:multiblock_graph} illustrates the construction.

\input{figures/figure_multiblock_graph}

We index intersection vertices by aisle \(j \in [1,m]\) and cross-aisle \(k \in [1,n]\), writing \(v_{j,k} \in V\).
The pick location set is \(P = \{p_0, p_1, \ldots, p_q\}\), where \(p_0\) is the depot and \(\{p_1,\ldots,p_q\}\) are item locations.
The depot may coincide with an intersection vertex; other pick locations are located within subaisles.
Edges \(E\) represent traversable aisle and cross-aisle segments, and their weights are the corresponding rectilinear distances.
Following standard practice for minimal-tour modelling \citep{ratliff1983order}, the graph is constructed to avoid redundant parallel edges.

\subsection{Tours and tour subgraphs}
\label{sec:tours}

The \SPRP{} seeks a minimum-length closed walk that starts and ends at the depot and visits every required pick location.
Equivalently, it seeks a minimum-weight \emph{tour subgraph} \(T \subseteq G\) from which such a tour can be constructed.
Figure~\ref{fig:multiblock_subgraph} shows an example tour subgraph.

\input{figures/figure_multiblock_subgraph}

Classical Eulerian graph theory implies that a subgraph \(T\) admits a closed walk traversing each of its edges exactly once
if and only if it is connected and all degrees are even.
Thus, a subgraph \(T \subseteq G\) is a tour subgraph for the \SPRP{} if and only if it (i) contains all required pick vertices, (ii) is connected, and (iii) has even degree at every vertex \citep{ratliff1983order}.

\paragraph{How we use structural insights}
Rectangular warehouse tours are not arbitrary: optimal solutions tend to follow simple geometric patterns that repeat across instances.
Rather than treating the warehouse as a generic graph and then relying on large families of subtour-elimination constraints, our formulations encode these patterns directly in the decision variables and connectivity constraints.
Concretely, we represent horizontal movement between adjacent aisles using a small set of cross-aisle configuration indicators and represent within-aisle movement using contiguous ``pass'' and ``branch-and-pick'' structures.
These design choices reduce the number of configurations the solver must consider and limit symmetry, which strengthens the relaxation and improves scalability in practice.
We reference prior analyses (including \citet{dunn2025double} and related manuscripts under review) only to provide intuition for why these choices are well-aligned with rectangular layouts; the mathematical programs in this paper are fully specified and can be checked for feasibility and correctness from the constraints provided here.


\subsection{Single Picker Routing Problem with Scattered Storage}

In the standard SPRP, each SKU is stored in a single location
within the warehouse,
with the picker only required to determine the best path
that visits all required pick positions.
Under \emph{scattered storage}, by contrast, one or more SKUs are available
at multiple positions throughout the warehouse \citep{weidinger2019picker}.
The picker must then decide not only how to route
but also where to pick each requested SKU when several positions offer it.
Scattered storage has become common in practice because this flexibility
often allows required items to be collected closer to the picker's starting position,
yielding shorter tours \cite{boysen2019warehousing},
but at the cost of greater routing complexity.

The resulting problem is referred to as
the single picker routing problem with scattered storage (SPRP-SS)
where the model must jointly decide which pick positions to use and how to route
a tour that visits them.
Figure~\ref{fig:scattered_tour} illustrates the setting
where some SKUs have multiple available positions in the warehouse,
but the tour need only visit one position per SKU (or enough to satisfy demand),
so not all positions holding a requested SKU are visited.

\input{figures/figure_scattered_tour}

The SPRP-SS has been proven to be NP-hard \citep{weidinger2018picker},
and remains so even when tours are restricted to
simple routing policies such as traversal,
return, midpoint, largest gap, and composite \citep{luke2024single}.
The problem can be further categorised as \\emph{unit demand},
where at most one unit per SKU is required,
or \\emph{general demand}, where a required quantity per SKU is specified.
This means that under general demand,
the picker may need to visit several positions
for a single SKU to satisfy its demand.
For unit-demand, the SPRP-SS is a special case
of the generalised travelling salesman problem (GTSP) \citep{daniels1998model},
but although solving via traditional GTSP formulations is convenient,
these methods do not exploit the structural properties of the warehouse layout.

We refer to Section~\ref{sec:related} for a brief overview of exact solution methods
and mathematical-programming formulations in the literature and how CC/EC are positioned among them.

%% file: figures/figure_multiblock_warehouse.tex
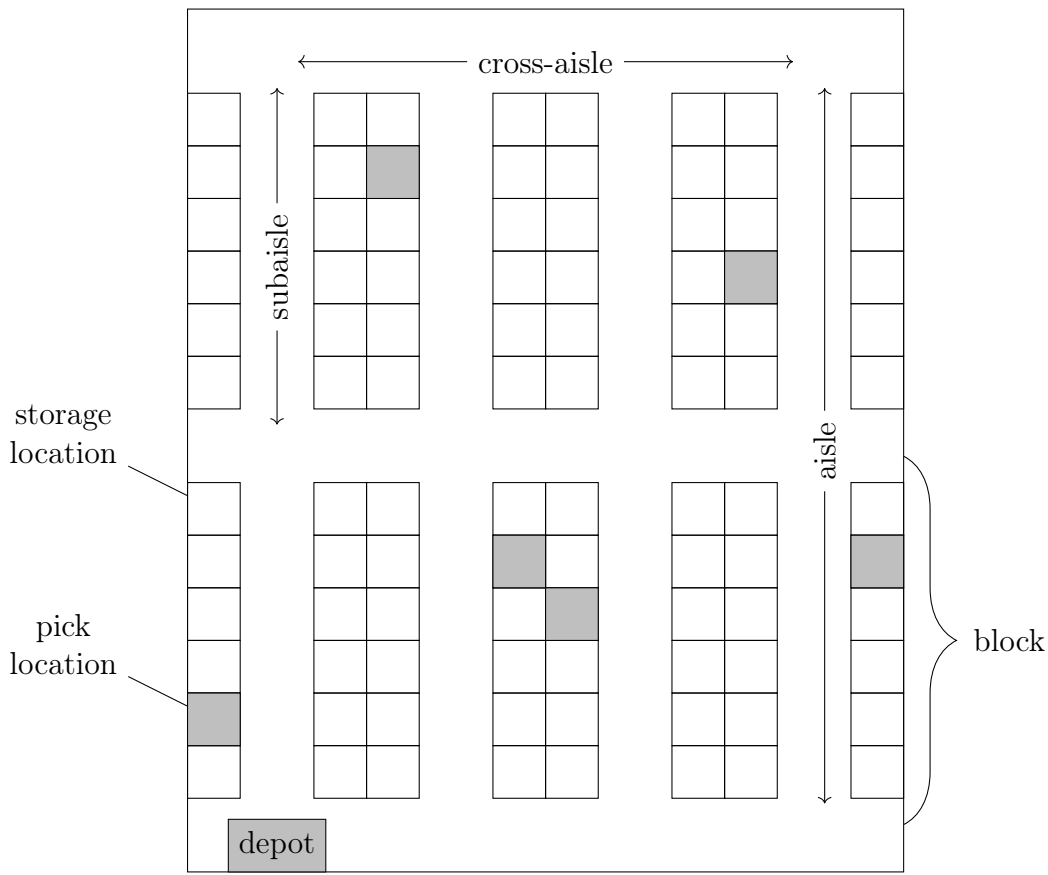
\begin{figure}[ht]
\centering

\begin{tikzpicture}

\tikzset{
  item/.style={rectangle, draw=black, fill=lightgray, minimum size=\nodesize cm
  }
}

\pgfmathsetmacro{\x}{0.7}
\pgfmathsetmacro{\o}{0.14}
\pgfmathsetmacro{\m}{4}
\pgfmathsetmacro{\n}{6}
\pgfmathsetmacro{\b}{2}

\pgfmathsetmacro{\nodesize}{0.7}

\foreach \j in {1,...,\m}{

    \foreach \k in {1,...,\b}{

        \foreach \i in {1,...,\n}{

            \node[rectangle,draw=black, minimum size=\nodesize cm] (l_\j\k\i) at
            (3*\x*\j+2*\o*\j-\x-\o, \i*\x+\k*\x*\n+\k*\x+2*\o*\k) { };

            \node[rectangle,draw=black, minimum size=\nodesize cm] (r_\j\k\i) at
            (3*\x*\j+2*\o*\j+\x+\o, \i*\x+\k*\x*\n+\k*\x+2*\o*\k) { };

        }

    }

}

\draw[black, thin]
(3*\x+2*\o-\x-\o-0.5*\x, \x+\x*\n+\x+2*\o-2*\o-1.5*\x) --
(3*\x*\m+2*\o*\m+\x+\o+0.5*\x, \x+\x*\n+\x+2*\o-2*\o-1.5*\x) --
(3*\x*\m+2*\o*\m+\x+\o+0.5*\x, \n*\x+\b*\x*\n+\b*\x+3*\o*\b+1.5*\x+\o) --
(3*\x+2*\o-\x-\o-0.5*\x, \n*\x+\b*\x*\n+\b*\x+3*\o*\b+1.5*\x+\o) -- cycle;

\node[item] (depot) at
(3*\x+2*\o, \x+\x*\n+\x+2*\o-2*\o-1.5*\x+0.5*\nodesize) {depot};

\node (cross)
at (1.5*\x*\m+\o*\m+\x+\o+0.5*\x, \n*\x+\b*\x*\n+\b*\x+3*\o*\b+0.5*\x+\o)
{cross-aisle};
\draw[->] (cross) -- (3*\x+4*\o, \n*\x+\b*\x*\n+\b*\x+3*\o*\b+0.5*\x+\o);
\draw[->] (cross) -- (3*\x*\m+2*\o*\m-2*\o, \n*\x+\b*\x*\n+\b*\x+3*\o*\b+0.5*\x+\o);

\node[rotate=90] (aisle)
at (3*\x*\m+\o*\m+\x, \n*\x+0.5*\b*\x*\n+0.5*\b*\x+2*\o*\b+\o+0.5*\x)
{aisle};
\draw[->] (aisle) -- (3*\x*\m+\o*\m+\x, 2*\x+\b*\x*\b*\x+3*\o*\b+\x+4*\o);
\draw[->] (aisle) -- (3*\x*\m+\o*\m+\x, \n*\x+\b*\x*\n+\b*\x+3*\o*\b+\o);

\node[rotate=90] (subaisle)
at (3*\x+2*\o, 0.5*\n*\x+\b*\x*\n+\b*\x+2*\o*\b+\o)
{subaisle};
\draw[->] (subaisle) -- (3*\x+2*\o, \n*\x+\b*\x*\n+\b*\x+3*\o*\b+\o);
\draw[->] (subaisle) -- (3*\x+2*\o, \n*\x+0.5*\b*\x*\n+0.5*\b*\x+2*\o*\b+\x+\o);

\draw [decorate, decoration={brace, mirror, amplitude=20pt}]
(3*\x*\m+2*\o*\m+\x+\o+0.5*\x, \x*\n+\x+2*\o) --
(3*\x*\m+2*\o*\m+\x+\o+0.5*\x, \x+\x*\n+\x+2*\o+\n*\x)
node[midway, xshift=40pt] {block};

\node[item] (item1) at (l_112) { };
\node[item] (item2) at (r_215) { };
\node[item] (item3) at (l_225) { };
\node[item] (item4) at (l_314) { };
\node[item] (item5) at (r_415) { };
\node[item] (item6) at (l_423) { };

\node[align=center] at ($(l_116)+(-2, 1)$) (storage) {storage\\location};
\draw[-] (storage) -- (l_116);

\node[align=center] at ($(item1)+(-2, 1)$) (item) {pick\\location};
\draw[-] (item) -- (item1);

\end{tikzpicture}
\caption{
    Rectangular warehouse layout with relevant features labelled.
}
\label{fig:multiblock_warehouse}
\end{figure}

%% file: figures/figure_multiblock_graph.tex
\begin{figure}[ht]
\centering

\begin{tikzpicture}

\tikzset{
  dbl/.style={-, thick, draw=black, double,
  double distance=8pt
  }
}

\tikzset{
  item/.style={circle, draw=black, fill=lightgray, minimum size=\nodesize pt
  }
}

\pgfmathsetmacro{\x}{0.75}
\pgfmathsetmacro{\m}{3}
\pgfmathsetmacro{\n}{6}
\pgfmathsetmacro{\b}{2}
\pgfmathsetmacro{\nodesize}{30}

\foreach \j in {0,...,\m}{

    \foreach \k in {0,...,\b}{

        \node[circle,draw=black,minimum size=\nodesize pt]
        (v-\j\k) at (2.5*\j*\x, \k*\n*\x) {\small $v_{\j,\k}$};

        \ifnum\k<\b
            \foreach \i in {1,...,\n}{
                \node (p_\j\k\i) at (2.5*\j*\x, \k*\n*\x+\i*\x) { };
            }
        \fi

    }

}

\node[item] (v-00) at (v-00) {\small $p_0$};

\node[item] (item1) at (p_002) {\small $p_1$};
\node[item] (item2) at (p_104) {\small $p_2$};
\node[item] (item3) at (p_114) {\small $p_3$};
\node[item] (item4) at (p_203) {\small $p_4$};
\node[item] (item5) at (p_304) {\small $p_5$};
\node[item] (item6) at (p_313) {\small $p_6$};

\draw[dbl] (v-00) -- (item1) -- (v-01) -- (v-02);
\draw[dbl] (v-10) -- (item2) -- (v-11) -- (item3) -- (v-12);
\draw[dbl] (v-20) -- (item4) -- (v-21) -- (v-22);
\draw[dbl] (v-30) -- (item5) -- (v-31) -- (item6) -- (v-32);

\draw[dbl] (v-00) -- (v-10) -- (v-20) -- (v-30);
\draw[dbl] (v-01) -- (v-11) -- (v-21) -- (v-31);
\draw[dbl] (v-02) -- (v-12) -- (v-22) -- (v-32);

\end{tikzpicture}
\caption{
    Graph representation of the warehouse in Figure~\ref{fig:multiblock_warehouse}.
}
\label{fig:multiblock_graph}
\end{figure}
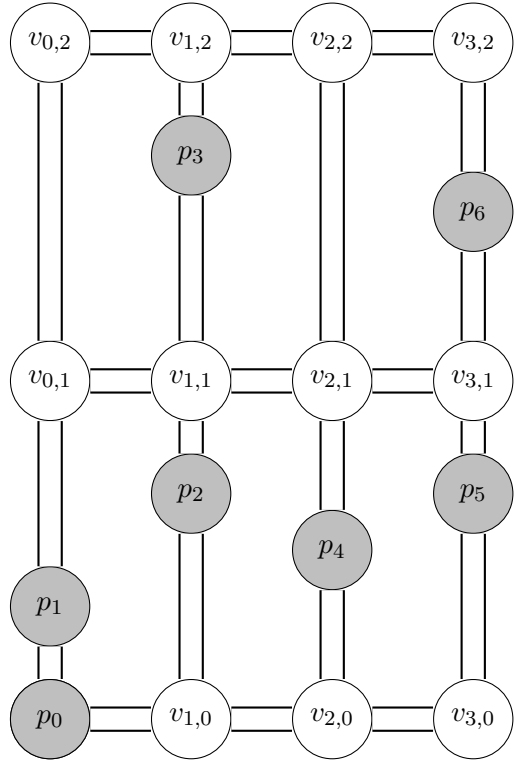

%% file: figures/figure_multiblock_subgraph.tex
\begin{figure}[ht]
\centering

\begin{tikzpicture}

\tikzset{
  dbl/.style={-, thick, draw=black, double,
  double distance=8pt
  }
}

\tikzset{
  sgl/.style={-, thick, draw=black
  }
}

\tikzset{
  item/.style={circle, draw=black, fill=lightgray, minimum size=\nodesize pt
  }
}

\pgfmathsetmacro{\x}{0.75}
\pgfmathsetmacro{\m}{3}
\pgfmathsetmacro{\n}{6}
\pgfmathsetmacro{\b}{2}
\pgfmathsetmacro{\nodesize}{30}

\foreach \j in {0,...,\m}{

    \foreach \k in {0,...,\b}{

        \node[circle,draw=black,minimum size=\nodesize pt]
        (v-\j\k) at (2.5*\j*\x, \k*\n*\x) {\small $v_{\j,\k}$};

        \ifnum\k<\b
            \foreach \i in {1,...,\n}{
                \node (p_\j\k\i) at (2.5*\j*\x, \k*\n*\x+\i*\x) { };
            }
        \fi

    }

}

\node[item] (v-00) at (v-00) {\small $p_0$};

\node[item] (item1) at (p_002) {\small $p_1$};
\node[item] (item2) at (p_104) {\small $p_2$};
\node[item] (item3) at (p_114) {\small $p_3$};
\node[item] (item4) at (p_203) {\small $p_4$};
\node[item] (item5) at (p_304) {\small $p_5$};
\node[item] (item6) at (p_313) {\small $p_6$};

\draw[dbl] (v-00) -- (item1);
\draw[sgl] (v-10) -- (item2) -- (v-11) -- (item3) -- (v-12);
\draw[sgl] (v-20) -- (item4) -- (v-21);
\draw[dbl] (item5) -- (v-31);
\draw[sgl] (v-31) -- (item6) -- (v-32);

\draw[dbl] (v-00) -- (v-10);
\draw[sgl] (v-10) -- (v-20);
\draw[sgl] (v-21) -- (v-31);
\draw[sgl] (v-12) -- (v-22) -- (v-32);

\end{tikzpicture}
\caption{
    An example of a tour subgraph for the warehouse in Figure~\ref{fig:multiblock_warehouse}.
}
\label{fig:multiblock_subgraph}
\end{figure}
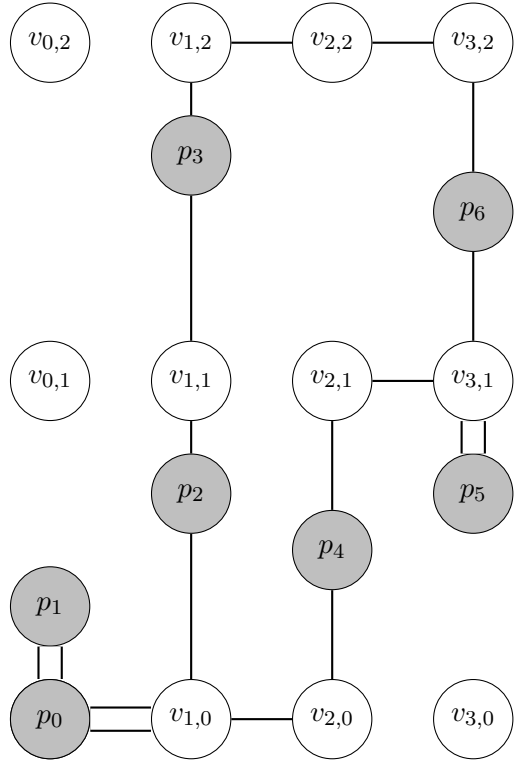

%% file: figures/figure_scattered_tour.tex
\begin{figure}[ht]
    \centering
    
    \begin{tikzpicture}
    
    \tikzset{
      dbl/.style={-, thick, draw=black, double,
      double distance=8pt
      }
    }
    
    \tikzset{
      sgl/.style={-, thick, draw=black
      }
    }
    
    \tikzset{
      item/.style={circle, draw=black, fill=lightgray, minimum size=\nodesize pt
      }
    }
    
    \pgfmathsetmacro{\x}{0.75}
    \pgfmathsetmacro{\m}{3}
    \pgfmathsetmacro{\n}{6} 
    \pgfmathsetmacro{\b}{2} 
    \pgfmathsetmacro{\nodesize}{30}
    
    \foreach \j in {0,...,\m}{
    
    
        \foreach \k in {0,...,\b}{    
    
            \node[circle,draw=black,minimum size=\nodesize pt]
            (v-\j\k) at (2.5*\j*\x, \k*\n*\x) {\small $v_{\j,\k}$};
    
            \ifnum\k<\b
                \foreach \i in {1,...,\n}{
                    \node (p_\j\k\i) at (2.5*\j*\x, \k*\n*\x+\i*\x) { };
                }
            \fi
    
        }
    
    }
    
    \node[item] (v-00) at (v-00) {\small $d$};
    
    \node[item] (item1) at (p_002) {\small A};
    \node[item] (item2) at (p_104) {\small B};
    \node[item] (item3) at (p_114) {\small A};
    \node[item] (item4) at (p_203) {\small B};
    \node[item] (item5) at (p_304) {\small A};
    \node[item] (item6) at (p_313) {\small C};
    
    \draw[dbl] (v-21) -- (v-31) -- (item6);
    \draw[sgl] (v-00) -- (item1) -- (v-01) -- (v-11) -- (v-21)
    -- (item4) -- (v-20) -- (v-10) -- (v-00);
    
    \end{tikzpicture}
    \caption[Example of a tour subgraph in a warehouse with scattered-storage.]
    {
        An example of a tour subgraph for a scattered-storage instance
        with SKUs A and B available at multiple positions.
    } 
    \label{fig:scattered_tour}
\end{figure}

%% file: sections_related.tex
This section positions our formulations within the order-picking literature, with emphasis on \emph{exact} methods for the \SPRP{} and its scattered-storage variant (\SPRPsS{}), and on why compact MILP formulations are useful for industrial decision-support settings.

\subsection{Exact methods for rectangular layouts}
Rectangular warehouses with parallel aisles have long supported specialised exact algorithms.
The dynamic-programming line initiated by \citet{ratliff1983order} exploits the layout geometry to solve the single-block \SPRP{}, and subsequent work extends these ideas to two-block and multi-block settings \citep{roodbergen2001routing,pansart2018exact}.
For larger layouts the problem becomes strongly $\mathcal{NP}$-hard \citep{prunet2025note}, motivating continued refinement of representations that retain tractability in practice \citep{hessler2022note}.
In this paper we do not develop new DP machinery; rather, we use DP-based approaches as performance benchmarks and as motivation for capturing layout structure directly within solver-friendly formulations.

\subsection{Mathematical-programming formulations}
Mathematical programming is particularly attractive when routing must be embedded in broader optimisation models (e.g., storage and slotting decisions, policy evaluation, or what-if layout analysis), and it becomes especially natural under scattered storage where pick-location choice and routing are coupled.
For single-block layouts, \citet{goeke2021modeling} propose an edge-configuration MILP (GS) that compactly represents admissible tour subgraphs and extends cleanly to \SPRPsS{}.
More recently, \citet{hessler2024exact} model a shortest path on a DP-derived state graph as a network-flow formulation (NF) and extend the state space to handle \SPRPsS{}, yielding a strong exact baseline for single- and two-block layouts.

Beyond rectangular single- and two-block settings, exact approaches include branch-and-cut schemes and formulations based on reductions to TSP variants and related routing problems \citep{valle2017optimally,saylam2024arc,irnich2025single,wildt2025picker}.
A recurring limitation of generic TSP encodings is that they treat the warehouse as an arbitrary graph; as a result, they typically require general subtour-elimination machinery and do not directly exploit the strong regularities induced by parallel-aisle geometry.

\subsection{Positioning of CC and EC (modelling contribution)}
The contribution of this paper is the development of \emph{structure-aware MILP formulations} for rectangular warehouses that reduce redundant configuration choices and simplify connectivity modelling.
Our \emph{Configuration Connectivity} (CC) formulation targets single-block layouts: it is derived from GS by removing configuration options that are unnecessary in rectangular layouts and by rewriting connectivity using a reduced set of configuration decisions, yielding a smaller and less symmetric MILP.
Our \emph{Edge Connectivity} (EC) formulation instead enforces connectivity directly on the warehouse graph using a unified core and layout-specific connection constraints, providing a single modelling framework that covers both single-block and two-block layouts and supports \SPRPsS{} in the same structure.

We compare CC and EC empirically with GS and NF.
The focus is on formulation size and solver performance (presolve reductions, root bounds, and runtimes), which are key practical criteria when the routing model is used as a reusable building block in industrial engineering applications rather than as a stand-alone algorithmic solver.

%% file: sections/section_cc.tex
This section presents the \emph{Configuration Connectivity} (CC) formulation, a compact MILP for the single-block \SPRP{} and \SPRPsS{} in rectangular warehouses.
CC is derived from the edge-configuration model of \citet{goeke2021modeling} (GS) but is intentionally \emph{smaller and less symmetric}: it removes configuration options that do not improve feasibility or cost in the layout class considered and it streamlines connectivity modelling by expressing it through a reduced set of horizontal configuration decisions.
The formulation is stated fully in this paper, so its feasibility semantics can be verified directly from the constraints; prior structural analyses \citep{dunn2025double,dunn2025deterministic} are cited only as background intuition for why the resulting model is compact and effective in rectangular layouts.
We refer to the model as \emph{Configuration Connectivity} because it enforces tour connectivity via cross-aisle configuration indicators in the spirit of the classical rectangular-layout structure \citep{ratliff1983order}.

The complete GS formulation is given in Appendix~\ref{sec:appendix_gs}.
We structure our CC formulation in the same way with
relevant model components appearing in the same order,
allowing us to clearly demonstrate the similarities and differences.

\subsection{Notation and Definitions}

\input{sections/section_cc_notation}

\subsection{Single-Block SPRP Formulation}

\input{sections/section_cc_sprp}

\subsection{Single-Block SPRP-SS Formulation}

\input{sections/section_cc_sprpss}

\subsection{Summary of Modifications}

The modifications above illustrate how the structural insights simplify the GS formulation.
Compared with GS, the CC formulation has fewer variables,
fewer constraints, and less redundancy among its constraints.
The key changes are as follows:

\begin{itemize}[noitemsep]
\item The double-traversal variable is removed entirely.
\item Two integer degree-parity variables per aisle are replaced by one binary variable.
\item Double traversals are removed from several constraints.
\item Some constraints are rewritten in terms of horizontal configurations only.
\item Two constraints are removed (one redundant, one merging into another).
\end{itemize}

\paragraph{Validity and modelling scope}
CC models the same rectangular single-block tour-subgraph requirements used in GS-style formulations: all required pick locations are covered, all relevant degrees satisfy even-parity conditions, and the selected edges/configurations form a single connected component.
The constraints in this section enforce these requirements using a reduced variable set, so feasibility and correctness can be assessed directly from the model as stated.
Appendix~\ref{sec:appendix_gs} records the GS baseline for comparison; the differences are in model size and constraint structure (and therefore in practical solvability), not in the underlying optimisation objective for the layout class considered.

%% file: sections/section_cc_notation.tex
The following notation follows \citet{goeke2021modeling} with
the four horizontal cross-aisle configurations
and the valid tour subgraph
requirements originating from \citet{ratliff1983order}.
The set $\Jcal = \{0, \ldots, m-1\}$ indexes $m$ aisles numbered from left to right.
Each aisle $j \in \Jcal$ has $n$ available picking positions numbered from bottom
to top and is associated with a set of required picking positions
$\Ical_j \subseteq \{0, \ldots, n-1\}$.
The depot is located at aisle $l$ with the binary variable $\thetad = 1$
if the depot is in the top cross-aisle and $0$ otherwise.

\paragraph{Cross-aisle configurations}
Four binary variables define the horizontal configuration between
aisles $j$ and $j+1$.
Specifically,
$\xx{x}_j$ = 1 if the bottom cross-aisle is traversed twice
(and top not used between $j$ and $j+1$),
$\xb{x}_j$ = 1 if the top cross-aisle
is traversed twice (bottom not used),
$\xbx{x}_j$ = 1 if both cross-aisles
are traversed once each and
$\xxb{x}_j$ = 1 if both cross-aisles
are traversed twice.

\paragraph{Vertical configurations}
Binary variables are also used to define vertical configurations with
$\xpass{j}$ representing aisle $j$ being traversed once ($1pass$) and
$\xppass{j}$ aisle $j$ traversed twice ($2pass$).
A significant difference between GS and CC is omitting the latter variable
to reduce redundancy and symmetry in rectangular layouts.
Variables $\pxji{x}{ji}$ and $\qxji{x}{ji}$ represent a vertical branch-and-pick from
the bottom and top cross-aisles, respectively, to position $i$ in aisle $j$.
Figure \ref{fig:formulations_gs} illustrates the vertical configuration variables.

\input{figures/figure_formulations_gs}

\paragraph{Degree parity}
GS uses integer variables $\pitop{j}$ and $\pibot{j}$ to represent the
total number of edges incident to
the top and bottom of aisle $j$, respectively (divided by two to ensure even parity).
CC replaces these with a single binary $\pi_j$ per aisle as we observe
that there are three variables that can introduce odd parity in aisle $j$,
specifically $\xpass{j}$, $\xbx{x}_{j-1}$, and $\xbx{x}_j$,
and enforcing a zero or even number of these variables is
sufficient to guarantee the degree parity condition is met.

\paragraph{Connectivity}
The $\tau_j$ binary variable in GS equals $1$ if the subgraph to the left
of aisle $j$ has two connected components.
This is used to ensure any separate components are propagated along the warehouse
to avoid creating disconnected subtours.
Once the last aisle is reached,
a single connected component is enforced.
Our CC formulation adopts this same logic
but is able to simplify and even remove some
constraints by expressing connectivity using a reduced set of horizontal configuration decisions that are sufficient for the rectangular layouts considered.

\paragraph{Costs}
Cost coefficients $c$ are precomputed for each configuration.
These represent the total cost of the edges introduced.
GS presents a unique cost for each configuration variable,
however, in rectangular warehouses the costs are uniform
so we factorize for clarity without affecting the objective.
We then have $\bar{c}_j$ for the horizontal distance from aisle $j$ to the next aisle
and $c^{\vert}$ for the vertical length of any aisle.

\paragraph{Other notation}
As in GS,
the notation $\cond{condition}{\cdots \quad}$ is used to represent
the bracketed terms
only being included when the condition is true.
This simplifies the expressions and makes the formulation clearer
in the cases where boundary constraints are not always present.

%% file: figures/figure_formulations_gs.tex
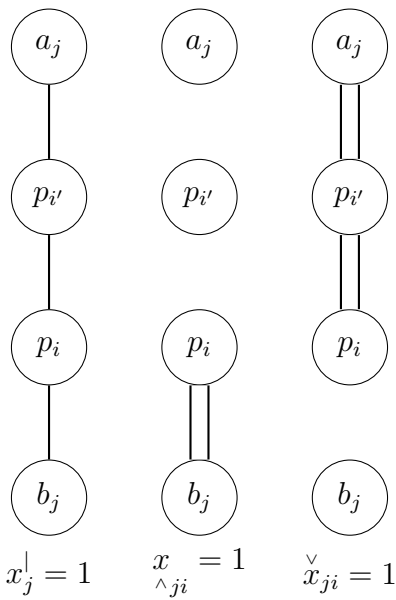
\begin{figure}[ht]
\centering
\begin{tikzpicture}[scale=1, every node/.style={circle, draw=black, minimum size=10mm}]

    \node (a1) at (0, 6) {$a_j$};
    \node (vi1) at (0, 2) {$p_i$};
    \node (vii1) at (0, 4) {$p_{i'}$};
    \node (b1) at (0, 0)   {$b_j$};
    \draw[thick] (a1) -- (vii1) -- (vi1) -- (b1);

    \node[draw=none] at (0, -1) {$\xpass{j}=1$};

    \node (a2) at (2, 6) {$a_j$};
    \node (vi2) at (2, 2) {$p_i$};
    \node (vii2) at (2, 4) {$p_{i'}$};
    \node (b2) at (2, 0)   {$b_j$};
    \draw[thick, double, double distance=6pt] (vi2) -- (b2);

    \node[draw=none] at (2, -1) {$\pxji{x}{ji}=1$};

    \node (a3) at (4, 6) {$a_j$};
    \node (vi3) at (4, 2) {$p_i$};
    \node (vii3) at (4, 4) {$p_{i'}$};
    \node (b3) at (4, 0)   {$b_j$};
    \draw[thick, double, double distance=6pt] (a3) -- (vii3) -- (vi3);

    \node[draw=none] at (4, -1) {$\qxji{x}{ji}=1$};

\end{tikzpicture}
\caption{Illustration of the vertical configuration variables in the GS formulation.}
\label{fig:formulations_gs}
\end{figure}

%% file: sections/section_cc_sprp.tex
\input{sections/section_cc_sprp_equations}

%% file: sections/section_cc_sprp_equations.tex
In this section we develop the Configuration Connectivity (CC) formulation for SPRP.
As it is based on the GS formulation of \citet{goeke2021modeling},
we will present it in the same structure
and address the changes to the GS formulation,
as well as acknowledge the components that
remain unchanged.

\CC

\paragraph{Objective (Minimize tour length)}
Minimize total travel cost.
\textbf{Modified.} The unnecessary $\xppass{j}$ term is removed.
Also, for rectangular warehouses,
horizontal costs are uniform, so we can factorize for clarity
without affecting the objective.
\begin{equation}
\min \sum_{j \in \Jcal} 2\bar{c}_j \bigl( \xx{x}_j + \xb{x}_j + \xbx{x}_j + 2\xxb{x}_j \bigr)
+ \sum_{j \in \Jcal} c^{\vert} \xpass{j}
+ \sum_{j \in \Jcal} \sum_{i \in \Ical_j}
\left( \pxji{c}{ji} \pxji{x}{ji} + \qxji{c}{ji} \qxji{x}{ji} \right)
\modeltag
\end{equation}

\paragraph{Cross-aisle configurations}
Ensure each aisle in the warehouse is visited using one of the four
cross-aisle configurations.
\textbf{Unchanged.}
\begin{equation}
\mbox{\begin{tabular*}{\dimexpr\linewidth-5em\relax}{@{\extracolsep{\fill}}lr@{}}%
$\displaystyle \xx{x}_j + \xb{x}_j + \xbx{x}_j + \xxb{x}_j = 1$ &
$\displaystyle j \in \Jcal \setminus \{m-1\}$%
\end{tabular*}}
\modeltag
\end{equation}

\paragraph{Visit all required positions}
Ensure the picker visits every required pick position
via a vertical configuration.
\textbf{Modified.} Variable $\xppass{j}$ has been removed as it is not required in a rectangular warehouse.
\begin{equation}
\mbox{\begin{tabular*}{\dimexpr\linewidth-5em\relax}{@{\extracolsep{\fill}}lr@{}}%
$\displaystyle \xpass{j} + \sum_{i' \in \Ical_j : i' \geq i} \qxji{x}{ji'} + \sum_{i' \in \Ical_j : i' \leq i} \pxji{x}{ji'} \geq 1$ &
$\displaystyle j \in \Jcal,\; i \in \Ical_j$%
\end{tabular*}}
\modeltag
\end{equation}

\paragraph{Top and bottom connectivity}
Allow a vertical branch-and-pick from the top or bottom cross-aisle
only if that cross-aisle is connected horizontally to an adjacent aisle.
This prevents isolated subtours within an aisle.
The depot is excluded from this constraint to allow for single-aisle tours,
with depot connectivity addressed exclusively later.
\textbf{Unchanged.}
\begin{equation}
\cond{j>0}{\xx{x}_{j-1} + \xbx{x}_{j-1} + \xxb{x}_{j-1}} + \xx{x}_j + \xbx{x}_j + \xxb{x}_j
\geq \pxji{x}{ji}
\qquad \begin{cases} \text{if } \thetad = 1: \; j \in \Jcal \\ \text{else: } j \in \Jcal \setminus \{l\} \end{cases},\; i \in \Ical_j
\modeltag
\end{equation}
\begin{equation}
\cond{j>0}{\xb{x}_{j-1} + \xbx{x}_{j-1} + \xxb{x}_{j-1}} + \xb{x}_j + \xbx{x}_j + \xxb{x}_j
\geq \qxji{x}{ji}
\qquad \begin{cases} \text{if } \thetad = 0: \; j \in \Jcal \\ \text{else: } j \in \Jcal \setminus \{l\} \end{cases},\; i \in \Ical_j
\modeltag
\end{equation}

\paragraph{Switches between top and bottom}
GS ensures switches between $\xx{x}$ and $\xb{x}$
are connected via $\xppass{j}$.
\textbf{Modified.}
With $\xppass{j}$ removed,
we can forbid such transitions entirely.
\begin{equation}
\mbox{\begin{tabular*}{\dimexpr\linewidth-5em\relax}{@{\extracolsep{\fill}}lr@{}}%
$\displaystyle \xx{x}_{j-1} + \xb{x}_j \leq 1$ &
$\displaystyle j \in \Jcal \setminus \{0\}$%
\end{tabular*}}
\modeltag
\end{equation}
\begin{equation}
\mbox{\begin{tabular*}{\dimexpr\linewidth-5em\relax}{@{\extracolsep{\fill}}lr@{}}%
$\displaystyle \xb{x}_{j-1} + \xx{x}_j \leq 1$ &
$\displaystyle j \in \Jcal \setminus \{0\}$%
\end{tabular*}}
\modeltag
\end{equation}

\paragraph{Depot inclusion}
GS ensures the depot is included in the tour via
a required cross-aisle configuration or traversals
$2\xppass{l}$ and $\xpass{l}$.
\textbf{Modified.}
With $\xppass{j}$ removed,
we also observe that $\xpass{l}$ can only be used in conjunction with
a $\xbx{x}_j$ cross-aisle configuration, which would connect the
depot regardless its location (top or bottom).
This allows us to simplify the depot inclusion constraints
to only consider cross-aisle configurations.
For the depot in the top cross-aisle ($\thetad = 1$).
\begin{equation}
\cond{l>0}{\xbx{x}_{l-1} + \xb{x}_{l-1} + \xxb{x}_{l-1}} + \xbx{x}_l + \xb{x}_l + \xxb{x}_l
\geq \cond{l>0}{\xx{x}_{l-1}} + \xx{x}_l
\modeltag
\end{equation}
For the depot in the bottom cross-aisle ($\thetad = 0$):
\begin{equation}
\cond{l>0}{\xbx{x}_{l-1} + \xx{x}_{l-1} + \xxb{x}_{l-1}} + \xbx{x}_l + \xx{x}_l + \xxb{x}_l
\geq \cond{l>0}{\xb{x}_{l-1}} + \xb{x}_l
\modeltag
\end{equation}

\paragraph{Degree parity}
GS uses $\pitop{j}$, $\pibot{j}$ to enforce even degree at each aisle warehouse vertex.
\textbf{Modified (and one constraint removed).}
A single binary $\pi_j$ suffices as we observe
that there are three variables that can introduce odd parity in aisle $j$,
$\xpass{j}$, $\xbx{x}_{j-1}$, and $\xxb{x}_j$.
Enforcing a zero or even number of these variables only requires a single binary variable.
This replaces the two individual equations for $\pitop{j}$ and $\pibot{j}$
with a single equation.
\begin{equation}
\mbox{\begin{tabular*}{\dimexpr\linewidth-5em\relax}{@{\extracolsep{\fill}}lr@{}}%
$\displaystyle \cond{j>0}{\xbx{x}_{j-1} + } \; \xbx{x}_j + \xpass{j} = 2\pi_j$ &
$\displaystyle j \in \Jcal$%
\end{tabular*}}
\modeltag
\end{equation}

\paragraph{Connected components}
Set $\tau_j = 1$ when two components arise
This occurs when $\xbx{x}_{j-1}$ or $\xxb{x}_{j-1}$
transitions to $\xxb{x}_j$ without connecting top and bottom in aisle $j$.
\textbf{Modified.}
Instead of considering both vertical and horizontal edges,
the deterministic structural result
shows that we can express connectivity by using horizontal variables only.
\begin{equation}
\mbox{\begin{tabular*}{\dimexpr\linewidth-5em\relax}{@{\extracolsep{\fill}}lr@{}}%
$\displaystyle \xxb{x}_j - \cond{j>0}{\xbx{x}_{j-1} + \xxb{x}_{j-1}} \leq \tau_j$ &
$\displaystyle j \in \Jcal$%
\end{tabular*}}
\modeltag
\end{equation}

\paragraph{Component propagation when left part unvisited}
Propagate the number of components along the warehouse.
\textbf{Modified.}
As above,
we simplify by propagating connectivity
with horizontal variables only.
\begin{equation}
\mbox{\begin{tabular*}{\dimexpr\linewidth-5em\relax}{@{\extracolsep{\fill}}lr@{}}%
$\displaystyle \tau_{j-1} - \xbx{x}_j \leq \tau_j$ &
$\displaystyle j \in \Jcal \setminus \{0\}$%
\end{tabular*}}
\modeltag
\end{equation}

\textbf{Propagate components}
GS required another constraint to propagate the number of components along the warehouse
that was dependent on the vertical configurations.
\textbf{Removed.}
That was made redundant by changing
the previous two constraints
to depend only on the horizontal variables.

\paragraph{Ensure $\xxb{x}_j$ when two components exist}
To prevent isolated subtours,
force $\xxb{x}_j = 1$ whenever there are two components.
\textbf{Unchanged.}
\begin{equation}
\mbox{\begin{tabular*}{\dimexpr\linewidth-5em\relax}{@{\extracolsep{\fill}}lr@{}}%
$\displaystyle \tau_j \leq \xxb{x}_j$ &
$\displaystyle j \in \Jcal$%
\end{tabular*}}
\modeltag
\end{equation}

\paragraph{Variable domains}
Define the decision variable domains.
\textbf{Modified.}
We removed $\xppass{j}$ as it is unnecessary and replace $\pitop{j}$, $\pibot{j}$,
two integer variables per aisle,
with a single binary variable $\pi_j$ to represent the degree parity.
\begin{equation}
\mbox{\begin{tabular*}{\dimexpr\linewidth-5em\relax}{@{\extracolsep{\fill}}lr@{}}%
$\displaystyle \xx{x}_j, \xb{x}_j, \xbx{x}_j, \xxb{x}_j, \tau_j \in \{0,1\}$ &
$\displaystyle j \in \Jcal \setminus \{m-1\}$%
\end{tabular*}}
\modeltag
\end{equation}
\begin{equation}
\mbox{\begin{tabular*}{\dimexpr\linewidth-5em\relax}{@{\extracolsep{\fill}}lr@{}}%
$\displaystyle \xpass{j} \in \{0,1\}$ &
$\displaystyle j \in \Jcal$%
\end{tabular*}}
\modeltag
\end{equation}
\begin{equation}
\mbox{\begin{tabular*}{\dimexpr\linewidth-5em\relax}{@{\extracolsep{\fill}}lr@{}}%
$\displaystyle \pxji{x}{ji}, \qxji{x}{ji} \in \{0,1\}$ &
$\displaystyle j \in \Jcal,\; i \in \Ical_j$%
\end{tabular*}}
\modeltag
\end{equation}
\begin{equation}
\mbox{\begin{tabular*}{\dimexpr\linewidth-5em\relax}{@{\extracolsep{\fill}}lr@{}}%
$\displaystyle \pi_j \in \{0,1\}$ &
$\displaystyle j \in \Jcal$%
\end{tabular*}}
\modeltag
\end{equation}

\paragraph{Last aisle}
The horizontal variables for the last aisle
are fixed to zero
as the warehouse is complete.
The connection variable is also zero for the last aisle
as the solution must contain a single connected component.
\textbf{Unchanged.}
\begin{equation}
\mbox{\begin{tabular*}{\dimexpr\linewidth-5em\relax}{@{\extracolsep{\fill}}l@{}}%
$\displaystyle \xx{x}_{m-1} = \xb{x}_{m-1} = \xbx{x}_{m-1} = \xxb{x}_{m-1} = \tau_{m-1} = 0$%
\end{tabular*}}
\modeltag
\end{equation}

%% file: sections/section_cc_sprpss.tex
\input{sections/section_cc_sprpss_equations}

%% file: sections/section_cc_sprpss_equations.tex
We now extend the CC formulation to the scattered-storage variant
using the same modelling approach as GS,
with the key modification that the double traversal ($2pass$) variable is excluded
in the visit-coverage constraints to keep the model compact and avoid redundant traversal patterns in rectangular layouts.
In \SPRPsS{}, a feasible plan jointly chooses which candidate storage positions supply each SKU and constructs a connected subtour through the visited positions subject to the same parity requirements as in the standard \SPRP{}; the constraints below link demand coverage, visit indicators, and the connectivity structure inherited from the formulation above.

For SPRP-SS, each SKU can be stored at multiple positions.
We use the same notation as \citet{goeke2021modeling} where
$\Hcal$ is the set of SKUs to pick,
$b_h$ is the demand for SKU $h \in \Hcal$,
$\Ical_{jh}$ contains picking positions in aisle $j$ from which SKU $h$ is available, and
$q_{jih}$ is the supply of SKU $h$ at position $i$ in aisle $j$.
The set $\Ical_j$ is redefined as all
positions in aisle $j$ from which any requested SKU is available,
not just the required pick positions.
Binary variable $x_{ji}$ indicates whether position $(j,i)$ is visited and
$\tilde{x}_j$ indicates whether aisle $j$ is reached by the picker.
The objective is the same as standard SPRP with $\Ical_j$ redefined as above.
The SPRP-SS formulation is then achieved by
replacing constraints (CC-2) and (CC-3) with the following.

\paragraph{Demand coverage}
Ensure the requested quantity of each SKU is picked from available positions.
When $b_h=1$ for every requested SKU and supply at each visited position is at most one unit, a single visit per SKU suffices; more generally, demand may be split across several positions.
\textbf{Unchanged.}
\begin{equation}
\mbox{\begin{tabular*}{\dimexpr\linewidth-5em\relax}{@{\extracolsep{\fill}}lr@{}}%
$\displaystyle \sum_{j \in \Jcal} \sum_{i \in \Ical_{jh}} q_{jih} x_{ji} \geq b_h$ &
$\displaystyle h \in \Hcal$%
\end{tabular*}}
\modeltag
\end{equation}

\paragraph{Visit selected positions}
Ensure every selected position $(j,i)$ is visited by the tour.
\textbf{Modified.} We remove $\xppass{j}$ as it is no longer required.
\begin{equation}
\mbox{\begin{tabular*}{\dimexpr\linewidth-5em\relax}{@{\extracolsep{\fill}}lr@{}}%
$\displaystyle \xpass{j} + \sum_{i' \in \Ical_j : i' \geq i} \qxji{x}{ji'} + \sum_{i' \in \Ical_j : i' \leq i} \pxji{x}{ji'} \geq x_{ji}$ &
$\displaystyle j \in \Jcal,\; i \in \Ical_j$%
\end{tabular*}}
\modeltag
\end{equation}

\paragraph{Position selection binary}
Define $x_{ji}$ as binary.
\textbf{Unchanged.}
\begin{equation}
\mbox{\begin{tabular*}{\dimexpr\linewidth-5em\relax}{@{\extracolsep{\fill}}lr@{}}%
$\displaystyle x_{ji} \in \{0,1\}$ &
$\displaystyle j \in \Jcal,\; i \in \Ical_j$%
\end{tabular*}}
\modeltag
\end{equation}

\paragraph{Aisle connectivity}
Ensure all aisles with selected positions (and the depot)
are visited.
\textbf{Unchanged.}
\begin{equation}
\mbox{\begin{tabular*}{\dimexpr\linewidth-5em\relax}{@{\extracolsep{\fill}}lr@{}}%
$\displaystyle \tilde{x}_j \geq x_{ji}$ &
$\displaystyle j \in \Jcal,\; i \in \Ical_j$%
\end{tabular*}}
\modeltag
\end{equation}
\begin{equation}
\mbox{\begin{tabular*}{\dimexpr\linewidth-5em\relax}{@{\extracolsep{\fill}}l@{}}%
$\displaystyle \tilde{x}_l = 1$%
\end{tabular*}}
\modeltag
\end{equation}
\begin{equation}
\mbox{\begin{tabular*}{\dimexpr\linewidth-5em\relax}{@{\extracolsep{\fill}}lr@{}}%
$\displaystyle \xx{x}_j + \xb{x}_j + \xbx{x}_j + \xxb{x}_j = \tilde{x}_{j+1}$ &
$\displaystyle j \in \Jcal \setminus \{m-1\} : j \geq l$%
\end{tabular*}}
\modeltag
\end{equation}
\begin{equation}
\mbox{\begin{tabular*}{\dimexpr\linewidth-5em\relax}{@{\extracolsep{\fill}}lr@{}}%
$\displaystyle \xx{x}_j + \xb{x}_j + \xbx{x}_j + \xxb{x}_j = \tilde{x}_j$ &
$\displaystyle j \in \Jcal : j < l$%
\end{tabular*}}
\modeltag
\end{equation}
\begin{equation}
\mbox{\begin{tabular*}{\dimexpr\linewidth-5em\relax}{@{\extracolsep{\fill}}lr@{}}%
$\displaystyle \tilde{x}_j \geq \tilde{x}_{j+1}$ &
$\displaystyle j \in \Jcal \setminus \{m-1\} : j \geq l$%
\end{tabular*}}
\modeltag
\end{equation}
\begin{equation}
\mbox{\begin{tabular*}{\dimexpr\linewidth-5em\relax}{@{\extracolsep{\fill}}lr@{}}%
$\displaystyle \tilde{x}_j \leq \tilde{x}_{j+1}$ &
$\displaystyle j \in \Jcal : j < l$%
\end{tabular*}}
\modeltag
\end{equation}
\begin{equation}
\mbox{\begin{tabular*}{\dimexpr\linewidth-5em\relax}{@{\extracolsep{\fill}}lr@{}}%
$\displaystyle \tilde{x}_j \in \{0,1\}$ &
$\displaystyle j \in \Jcal$%
\end{tabular*}}
\modeltag
\end{equation}

%% file: sections/section_ec.tex
We now present the \emph{Edge Connectivity} (EC) formulation, a graph-based MILP designed to handle single-block and two-block rectangular warehouse layouts within a unified modelling framework.
The formulation is inspired by the GS and CC modelling style but is organised around a \emph{common core} (objective, demand coverage, and degree-parity structure) together with \emph{layout-specific connection constraints} that enforce connectivity of the tour subgraph (single-block versus two-block).

The key modelling idea is to construct a tour by selecting a subgraph that (i) satisfies even-degree parity, (ii) covers all required picking positions, and (iii) is connected.
To keep the formulation compact and solver-friendly, EC encodes within-aisle movement through single-traversal and branch-and-pick style vertical configurations together with the horizontal edge structure.
These choices are well-aligned with the regularities observed in optimal tours for rectangular layouts; we cite \citet{dunn2025double} and related manuscripts under review for background intuition, but the formulation below is fully specified and its feasibility semantics follow from the constraints stated here.
Unlike CC, which is restricted to single-block layouts and encodes connectivity implicitly via configuration variables, EC enforces connectivity directly on the warehouse graph and only specialises the connection constraints by block type, yielding a single formulation family for single- and two-block layouts (and a natural pathway to additional blocks).

\paragraph{Modelling scope}
EC models tour subgraphs in the rectilinear warehouse graph: all required pick vertices are included, all relevant degrees satisfy even parity, and the selected edges form one connected component.
The connection constraints in the sections below are tailored to single-block and two-block layouts; the contribution is the resulting unified MILP structure and its computational performance relative to GS/CC-style MILPs and the NF state-graph benchmark.

\subsection{Considerations for General Formulation}
\input{sections/section_ec_prelim}
\subsection{Core Formulation}
\input{sections/section_ec_core}

\subsection{Connectivity Structure}

We now present the common connectivity components before
specifically addressing the constraints for single-block,
two-block layouts, and
a roadmap to applications in larger warehouses.

\paragraph{Connectivity}
In line with GS and CC, we maintain $\tau$ for connectivity but extend it
in EC to $\acr{j}{k}{k'}$ indicating whether crosses $k$ and $k'$ are connected at aisle $j$,
with (aisle, cross-aisle, cross-aisle) indices instead of the single aisle index $\tau_j$
in the single-block formulations.
\begin{equation}
\mbox{\begin{tabular*}{\dimexpr\linewidth-5em\relax}{@{\extracolsep{\fill}}lr@{}}%
$\displaystyle \acr{j}{k}{k'} \in [0,1]$ &
$\displaystyle j \in \Jcal,\; (k,k') \in \mathcal{P}; k < k'$%
\end{tabular*}}
\modeltag
\end{equation}

\paragraph{Previous-aisle connection}
The connection constraints for single-block and two-block layouts will be covered
in the next section;
however, we will introduce the previous-aisle connection here
as it is a common constraint in both single-block and two-block.
If $\acp{j}{k}{k'} = 1$ then at aisle $j-1$ both $k$
and $k'$ have at least one horizontal edge and $(k,k')$ is connected there.
\begin{equation}
\mbox{\begin{tabular*}{\dimexpr\linewidth-5em\relax}{@{\extracolsep{\fill}}lr@{}}%
$\displaystyle \acp{j}{k}{k'} \le \acbarx{j-1}{k} + \acbarbarx{j-1}{k}$ &
$\displaystyle j \in \Jcal \setminus \{0\},\; (k,k') \in \mathcal{P}$%
\end{tabular*}}
\modeltag
\end{equation}
\begin{equation}
\mbox{\begin{tabular*}{\dimexpr\linewidth-5em\relax}{@{\extracolsep{\fill}}lr@{}}%
$\displaystyle \acp{j}{k}{k'} \le \acbarx{j-1}{k'} + \acbarbarx{j-1}{k'}$ &
$\displaystyle j \in \Jcal \setminus \{0\},\; (k,k') \in \mathcal{P}$%
\end{tabular*}}
\modeltag
\end{equation}
\begin{equation}
\mbox{\begin{tabular*}{\dimexpr\linewidth-5em\relax}{@{\extracolsep{\fill}}lr@{}}%
$\displaystyle \acp{j}{k}{k'} \le \acr{j-1}{k}{k'}$ &
$\displaystyle j \in \Jcal \setminus \{0\},\; (k,k') \in \mathcal{P}$%
\end{tabular*}}
\modeltag
\end{equation}
\begin{equation}
\mbox{\begin{tabular*}{\dimexpr\linewidth-5em\relax}{@{\extracolsep{\fill}}lr@{}}%
$\displaystyle \acp{j}{k}{k'} \in [0,1]$ &
$\displaystyle j \in \Jcal,\; (k,k') \in \mathcal{P}; k < k'$%
\end{tabular*}}
\modeltag
\end{equation}

\paragraph{Last-aisle connection}
If a pair of cross-aisles in the last aisle both have a left edge,
then they must be connected to ensure a single connected tour subgraph.
\begin{equation}
\acbarx{m-2}{k} + \acbarbarx{m-2}{k} + \acbarx{m-2}{k'} + \acbarbarx{m-2}{k'} - 1
\le \acr{m-1}{k}{k'} \qquad (k,k') \in \mathcal{P}
\modeltag
\end{equation}

Integrality of connectivity variables
$\acr{j}{k}{k'}$ and  $\acp{j}{k}{k'}$ (as well as $\acz{j}{k}{k'}$ for two-block)
is implied by the rest of the formulation but not required for correctness.
These do not appear in the objective and
are constrained only by linear inequalities whose right-hand sides are sums.
They are therefore bounded (e.g.\ in $[0,1]$) and,
at any feasible solution consistent with the binary horizontal
and vertical decisions, the connectivity constraints force them
to take values 0 or 1.
Integrality of these variables is thus implied by the rest of the formulation,
so they can be declared continuous without changing the feasible set or the optimum.
Using continuous rather than binary variables for connectivity reduces the number
of integer variables in the MILP, which typically improves solution times as
the linear relaxations are unchanged,
but branch-and-bound explores fewer fractional combinations
and the solver can often prove optimality with less branching.

\subsubsection{Single-Block Connectivity}
\input{sections/section_ec_connectivity_singleblock}

\subsubsection{Two-Block Connectivity}
\input{sections/section_ec_connectivity_twoblock}

\subsection{Incorporating Scattered Storage}
\input{sections/section_ec_sprpss}

%% file: sections/section_ec_prelim.tex
The fundementals of this formulation,
including notation,
are kept as close to CC as possible,
however, some changes
are necessary to accommodate the generalisation to two-block layouts and beyond.
To allow for multiple blocks,
each block $k \in \Kcal \setminus \{ n-1 \}$ of aisle $j \in \Jcal$
is associated with a set of required picking positions
$\Ical_{jk} \subseteq \{k, \ldots, k+1\}$.

\paragraph{Cross-aisle configurations}
As opposed to GS and CC where complete horizontal configurations are encoded in a single variable,
EC encodes each cross-aisle configuration as a separately.
This is necessary to accommodate the generalisation to two-block layouts and beyond.
The variables are $\acbarx{j}{k}$ and $\acbarbarx{j}{k}$
for the single-edge and double-edge configurations in cross-aisle $k$ at aisle $j$, respectively.

\paragraph{Vertical configurations}
The vertical variables were also kept as close to CC as possible
with the necessary extentions for multiple blocks.
The variable $\acpass{j}{k}$ represents a full pass from cross $k$ to $k+1$
at aisle $j$, and is therefore not required in the last cross-aisle.
Variables $\qxji{x}{ji}$ and $\pxji{x}{ji}$ again represents a vertical pick
from the above and below, respectively,
to position $i$ in aisle $j$.
These have been modified though to only include the segments to
the next item or cross-aisle, whichever is closer.
This makes demand coverage more clear for each item as
we only need to check the specific variables associated with a location,
we do not need to look
for vertical configurations in items above or below.
It also means that for branch-and-pick connectivity,
we only need to check the item immediately above (or below) the cross-aisle,
rather than all items in the block.
Connectivity within the block is enforced by chain propagation
that ensures that if an item is picked from below,
so must all items below it until a cross-aisle is reached.
Figure \ref{fig:formulations_ec} illustrates the vertical configurations
and the concept of chain propagation.

\input{figures/figure_formulations_ec}

\paragraph{Depot}
Instead of $\thetad$ as a binary representation of the depot
being top or bottom,
we will let $\theta$ represent the location of the depot cross-aisle,
still assumed to be the top or bottom.
The full location is therefore given by the cross $\theta$ in aisle $l$.

\paragraph{Double traversals (more than one blocks)}
In the rectangular layouts studied here, allowing edges to be traversed twice can be useful as a modelling device, but it often introduces redundant configurations and symmetry that enlarge the MILP search space.
Accordingly, EC is designed around a compact set of single-traversal and branch-and-pick style components, with double traversals treated only as a special-case mechanism when they are unavoidable for depot access in very small layout regimes.
We cite the published elimination insights in \citet{dunn2025double} and related refinements (manuscripts under review) as background intuition for these restrictions; the formulation below is stated explicitly and should be read as a modelling choice that yields a compact solver-ready representation of the tour-subgraph requirements.
The remaining necessary double traversal still need to
be considered though,
which in two-blocks is only the case when all required (or selected)
pick positions are in the same aisle as the depot.
For the standard SPRP,
instances with one aisle can be can be
detected and solved trivially,
applying the formulation to those with more than one aisle only.
In the scattered-storage setting,
this is not so simple as SKUs can be stored in multiple aisles,
but the minimal solution may still involve a single aisle.
Figure \ref{fig:formulations_depot} illustrates
a necessary double traversal in a two-block SPRP as
well as two tours for the same SPRP-SS instance,
one containing multiple aisles and one containing only one aisle.
To this end, we extend the set of pick positions
to include the middle cross-aisle in the depot aisle,
$\Ical_l' = \Ical_l \cup \{v_{l,\theta}\}$,
letting us model the double traversal by
introducing the segment between the middle-cross and the
next pick location in the depot aisle.

\input{figures/figure_formulations_depot}

%% file: figures/figure_formulations_ec.tex
\begin{figure}[ht]
    \centering
    \begin{tikzpicture}[scale=1, every node/.style={circle, draw=black, minimum size=14mm}]
    
        \node (a1) at (0, 7.5) { };
        \node at (a1) {$v_{j,k+1}$};
        \node (vi1) at (0, 2.5) {$p_i$};
        \node (vii1) at (0, 5) {$p_{i'}$};
        \node (b1) at (0, 0)   { };
        \node at (b1) {$v_{j,k}$};
        \draw[thick] (a1) -- (vii1) -- (vi1) -- (b1);
    
        \node[draw=none] at (0, -1.5) {$\acpass{{j}}{{k}}=1$};
    
        \node (a2) at (3, 7.5) { };
        \node at (a2) {$v_{j,k+1}$};
        \node (vi2) at (3, 2.5) {$p_i$};
        \node (vii2) at (3, 5) {$p_{i'}$};
        \node (b2) at (3, 0)   { };
        \node at (b2) {$v_{j,k}$};
        \draw[thick, double, double distance=8pt] (b2) -- (vi2);
    
        \node[draw=none] at (3, -1.5) {$\pxji{{x}}{{ji}}=1$};
    
        \node (a3) at (6, 7.5) { };
        \node at (a3) {$v_{j,k+1}$};
        \node (vi3) at (6, 2.5) {$p_i$};
        \node (vii3) at (6, 5) {$p_{i'}$};
        \node (b3) at (6, 0) { };
        \node at (b3) {$v_{j,k}$};
        \draw[thick, double, double distance=8pt] (vii3) -- (vi3);
        \draw[color=red,thick, double, double distance=8pt] (a3) -- (vii3);
    
        \node[draw=none] at (6, -1.5) {$\qxji{{x}}{{ji}}=1$};
    
    \end{tikzpicture}
    \caption[Illustration of the vertical configuration variables in the EC formulation.]
    {
    Vertical configuration variables in the EC formulation for general warehouse layouts
    as shown by black lines.
    $\qxji{{x}}{{ji}}=1$ only includes the segment from $p_{i}$ to $p_{i'}$.
    Connectivity to the next cross-aisle as show by the red lines,
    is enforced by the chain propagation constraints.
    }
    \label{fig:formulations_ec}
\end{figure}
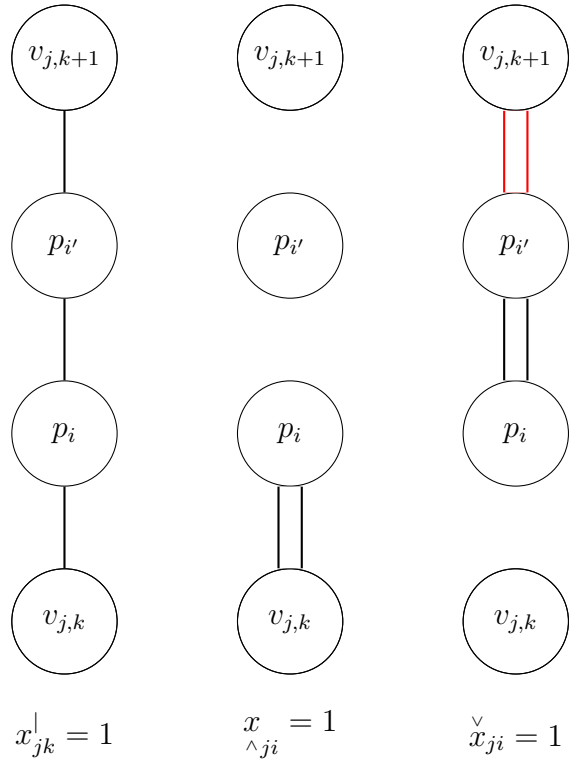

%% file: figures/figure_formulations_depot.tex
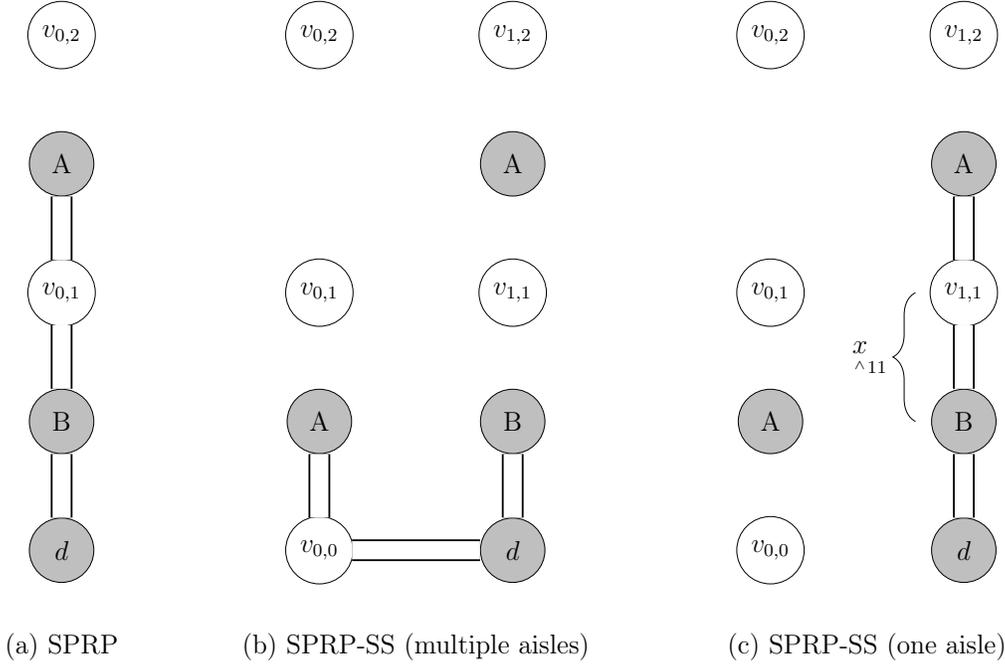
\begin{figure}[ht]
    \centering
    \resizebox{\linewidth}{!}{%
    \begin{tikzpicture}[
        scale=1,
        every node/.style={circle, draw=black, minimum size=10mm},
        execute at end picture={\path[use as bounding box] (-1, -2) rectangle (15, 8.7);}]
    \clip (-1, -2) rectangle (15, 8.7);
    
        \node[draw=none] (a1) at (0, 8) { };
        \node at (a1) {$v_{0,2}$};
        \node[fill=lightgray] (vi1) at (0, 6) {A};
        \node[fill=lightgray] (vii1) at (0, 2) {B};
        \node[draw=none] (b1) at (0, 4)   { };
        \node at (b1) {$v_{0,1}$};
        \node[draw=none] (c1) at (0, 0)   { };
        \node[fill=lightgray] at (c1) {$d$};
        \draw[thick, double, double distance=8pt] (vi1) -- (b1) -- (vii1)  -- (c1);
    
        \node[draw=none] at (0, -1.5) {(a) SPRP};
    
        \node[draw=none] (a2) at (4, 8) { };
        \node at (a2) {$v_{0,2}$};
        \node[fill=lightgray] (vii2) at (4, 2) {A};
        \node[draw=none] (b2) at (4, 4)   { };
        \node at (b2) {$v_{0,1}$};
        \node[draw=none] (c2) at (4, 0)   { };
        \node at (c2) {$v_{0,0}$};
        \draw[thick, double, double distance=8pt] (vii2)  -- (c2);
    
        \node[draw=none] (a3) at (7, 8) { };
        \node at (a3) {$v_{1,2}$};
        \node[fill=lightgray] (vi3) at (7, 6) {A};
        \node[fill=lightgray] (vii3) at (7, 2) {B};
        \node[draw=none] (b3) at (7, 4)   { };
        \node at (b3) {$v_{1,1}$};
        \node[draw=none] (c3) at (7, 0)   { };
        \node[fill=lightgray] at (c3) {$d$};
        \draw[thick, double, double distance=8pt] (vii3) -- (c3);

        \draw[thick, double, double distance=8pt] (c2) -- (c3);

        \node[draw=none] at (5.5, -1.5) {(b) SPRP-SS (multiple aisles)};

        \node[draw=none] (a4) at (11, 8) { };
        \node at (a4) {$v_{0,2}$};
        \node[fill=lightgray] (vii4) at (11, 2) {A};
        \node[draw=none] (b4) at (11, 4)   { };
        \node at (b4) {$v_{0,1}$};
        \node[draw=none] (c4) at (11, 0)   { };
        \node at (c4) {$v_{0,0}$};

        \node[draw=none] (a5) at (14, 8) { };
        \node at (a5) {$v_{1,2}$};
        \node[fill=lightgray] (vi5) at (14, 6) {A};
        \node[fill=lightgray] (vii5) at (14, 2) {B};
        \node[draw=none] (b5) at (14, 4)   { };
        \node at (b5) {$v_{1,1}$};
        \node[draw=none] (c5) at (14, 0)   { };
        \node[fill=lightgray] at (c5) {$d$};
        \draw[thick, double, double distance=8pt] (vi5) -- (b5) -- (vii5)  -- (c5);

        \node[draw=none] at (12.5, -1.5) {(c) SPRP-SS (one aisle)};

        \draw [decorate, decoration={brace, amplitude=10pt}]
        (13.25, 2) --
        (13.25, 4)
        node[draw=none, midway, xshift=-20pt] {$\pxji{x}{11}$};
    \end{tikzpicture}%
    }
    \caption[Examples of double traversals in scattered-storage instances.]
    {
    Example of a necessary double traversal in a two-block SPRP (a),
    a two-block SPRP-SS with a tour containing multiple aisles (b),
    and the same two-block SPRP-SS with a tour containing only one aisle (c).
    }
    \label{fig:formulations_depot}
\end{figure}

%% file: sections/section_ec_core.tex
The core EC formulation seeks a minimal cost tour while
enforcing four structural properties:
(i) visit of all required pick positions;
(ii) even degree parity at each vertex;
(iii) inner-aisle connectivity so that within each aisle the tour is connected;
and (iv) inter-aisle connectivity via the cross-aisle edges and the connection variables
that link crosses across aisles (so that the tour is a single connected component).
The equations below formalise these requirements.

\EC

\paragraph{Objective (Minimize tour length)}
Minimize total travel cost
of horizontal edges over all aisles and crosses,
single traversals over all aisles and blocks,
and vertical segment costs for all required positions
(including the depot middle-cross for two-blocks).
\begin{equation}
\min \sum_{j \in \Jcal} \sum_{k \in \Kcal} \bar{c}_j
\bigl( \acbarx{j}{k} + 2\acbarbarx{j}{k} \bigr)
+ \sum_{j \in \Jcal} \sum_{k \in \Kcal \setminus \{n-1\}} c^{\vert} \acpass{j}{k}
+ \sum_{j \in \Jcal} \sum_{i \in \Ical_j'}
\bigl( \pxji{c}{ji} \pxji{x}{ji} + \qxji{c}{ji} \qxji{x}{ji} \bigr)
\modeltag
\end{equation}

\paragraph{Cross-aisle configurations (optional)}
At most one of single-edge or double-edge variables can be active per cross.
This constraint is not required for a valid tour subgraph
as it is already enforced by degree parity and the minimal
objective function,
however it does provide tighter bounds on the variables.
\begin{equation}
\mbox{\begin{tabular*}{\dimexpr\linewidth-5em\relax}{@{\extracolsep{\fill}}lr@{}}%
$\displaystyle \acbarx{j}{k} + \acbarbarx{j}{k} \le 1$ &
$\displaystyle j \in \Jcal \setminus \{m-1\},\; k \in \Kcal$%
\end{tabular*}}
\modeltag
\end{equation}

\paragraph{Visit all required positions}
For the standard SPRP,
every required position $(j,i)$ is visited by
a single traversal in that subaisle,
a segment-from-below
or a segment-from-above.
The introduction of top and bottom segments allows
us to force this to exactly one instead
of having to check the variables of the pick positions
above and below.
Note that this applies to the original pick position set $\Ical_{jk}$
without the depot aisle middle-cross.
\begin{equation}
\mbox{\begin{tabular*}{\dimexpr\linewidth-5em\relax}{@{\extracolsep{\fill}}lr@{}}%
$\displaystyle \acpass{j}{c} + \pxji{x}{ji} + \qxji{x}{ji} = 1$ &
$\displaystyle j \in \Jcal,\; i \in \Ical_j$%
\end{tabular*}}
\modeltag
\end{equation}

\paragraph{Connectivity within a subaisle (Chain propagation)}
Connectivity within a subaisle requires a continuous chain of segments.
If a pick position is visited from below,
the previous pick position in the same block must also be visited from below.
\begin{equation}
\mbox{\begin{tabular*}{\dimexpr\linewidth-5em\relax}{@{\extracolsep{\fill}}lr@{}}%
$\displaystyle \pxji{x}{j(i+1)} \le \pxji{x}{ji}$ &
$\displaystyle j \in \Jcal,\; i, (i+1) \in \Ical_{jk},\; k \in \Kcal \setminus \{n-1\}$%
\end{tabular*}}
\modeltag
\end{equation}
If a pick position is visited from above,
the next pick position in the same block must also be visited from above.
\begin{equation}
\mbox{\begin{tabular*}{\dimexpr\linewidth-5em\relax}{@{\extracolsep{\fill}}lr@{}}%
$\displaystyle \qxji{x}{j(i-1)} \le \qxji{x}{ji}$ &
$\displaystyle j \in \Jcal,\; i, (i-1) \in \Ical_{jk},\; k \in \Kcal \setminus \{n-1\}$%
\end{tabular*}}
\modeltag
\end{equation}

\paragraph{Top and bottom connectivity}
Allow a vertical branch-and-pick from the top or bottom of a subaisle
only if that cross-aisle is connected horizontally to adjacent aisles.
As opposed to GS and CC,
who need to apply these constraints for each item in the subaisle,
we only need to consider the pick positions immediately above and below the
relevant cross-aisle.
\begin{equation}
\cond{j>0}{\acbarx{(j-1)}{k} + \acbarbarx{(j-1)}{k}}
+ \acbarx{j}{k} + \acbarbarx{j}{k} \geq \pxji{x}{ji}
\qquad
(j,k) \in \Jcal \times \Kcal \setminus \{(l,\theta)\},\; i = \min (\Ical_{jk}) 
\modeltag
\end{equation}
\begin{equation}
\cond{j>0}{\acbarx{(j-1)}{k} + \acbarbarx{(j-1)}{k}}
+ \acbarx{j}{k} + \acbarbarx{j}{k} \geq \qxji{x}{ji}
\qquad
(j,k) \in \Jcal \times \Kcal \setminus \{(l,\theta)\},\; i = \max (\Ical_{j(k-1)}) 
\modeltag
\end{equation}

\paragraph{Even number of horizontal edges (optional)}
We introduce an integer variable $\aceta{j}$
to enforce an even number of horizontal edges between aisles $j$ and $j+1$.
In the standard SPRP,
we can restrict this to $\aceta{j}>1$ as each aisle must have at least one horizontal edge pair.
Once again, this constraint is not required for a valid tour subgraph
as it is already enforced by degree parity and the upcomming
connectivity constraints, however it does enforce tighter bounds on the variables.
\begin{equation}
\mbox{\begin{tabular*}{\dimexpr\linewidth-5em\relax}{@{\extracolsep{\fill}}lr@{}}%
$\displaystyle \sum_{k \in \Kcal} (\acbarx{j}{k} + 2\acbarbarx{j}{k}) = 2\aceta{j}$ &
$\displaystyle j \in \Jcal \setminus \{m-1\}$%
\end{tabular*}}
\modeltag
\end{equation}

\paragraph{Depot inclusion}
The depot must be included in the tour
by ensuring it is connected by horizontal edges.
This is only enforced if there are other
horizontal edges to or from the depot aisle,
allowing for a single-aisle tour.
\begin{equation}
\cond{l>0}{\acbarx{(l-1)}{k} + \acbarbarx{(l-1)}{k}}
+ \acbarx{l}{k} + \acbarbarx{l}{k}
\ge
\cond{l>0}{\acbarx{(l-1)}{\theta} + \acbarbarx{(l-1)}{\theta}}
+\acbarx{l}{\theta} + \acbarbarx{l}{\theta}
\qquad k \in \Kcal \setminus \{\theta\}
\modeltag
\end{equation}

\paragraph{Degree parity}
At each cross vertex the degree is even.
Similar to CC,
we only need to consider the configurations
that introduce an odd number of edges.
The variable $\acpi{j}{k}$ is an integer 
representing half the number of single edge variables incident to $v_{j,k}$.
\begin{equation}
\acbarx{j}{k} \; \cond{j>0}{+ \acbarx{(j-1)}{k}} \; \cond{k>0}{+ \acpass{j}{(k-1)}} \cond{k<n-1}{+ \acpass{j}{k}}
= 2\acpi{j}{k}
\qquad j \in \Jcal,\; k \in \Kcal
\modeltag
\end{equation}

\paragraph{Variable domains}
All are binary except for degree parity and the optional edge-pair variable.
\begin{equation}
\mbox{\begin{tabular*}{\dimexpr\linewidth-5em\relax}{@{\extracolsep{\fill}}lr@{}}%
$\displaystyle \acbarx{j}{k}, \acbarbarx{j}{k}, \acpass{j}{k}, \qxji{x}{ji}, \pxji{x}{ji} \in \{0,1\}$ &
$\displaystyle j \in \Jcal,\; k \in \Kcal,\; i \in \Ical_j$%
\end{tabular*}}
\modeltag
\end{equation}
\begin{equation}
\mbox{\begin{tabular*}{\dimexpr\linewidth-5em\relax}{@{\extracolsep{\fill}}lr@{}}%
$\displaystyle \acpi{j}{k} \in \{0,2\}$ &
$\displaystyle j \in \Jcal,\; k \in \Kcal,\; i \in \Ical_j$%
\end{tabular*}}
\modeltag
\end{equation}
\begin{equation}
\mbox{\begin{tabular*}{\dimexpr\linewidth-5em\relax}{@{\extracolsep{\fill}}lr@{}}%
$\displaystyle \aceta{j} \in \{1, \ldots, |\Kcal|/2\}$ &
$\displaystyle j \in \Jcal \setminus \{m-1\}$%
\end{tabular*}}
\modeltag
\end{equation}

\paragraph{Last aisle}
The horizontal configurations in the last aisle are set to zero.
\begin{equation}
\mbox{\begin{tabular*}{\dimexpr\linewidth-5em\relax}{@{\extracolsep{\fill}}lr@{}}%
$\displaystyle \acbarx{m-1}{k} = \acbarbarx{m-1}{k} = 0$ &
$\displaystyle k \in \Kcal$%
\end{tabular*}}
\modeltag
\end{equation}

%% file: sections/section_ec_connectivity_singleblock.tex
The general formulation above
ensures a subgraph that has even degree parity and visits all required positions.
The last remaining requirement for a valid tour subgraph
is to ensure that the subgraph is connected.
We do this in a similar way to \citet{ratliff1983order} in that
we propagate connectivity from left to right.
As $2pass$ is not required for connectivity,
$1pass$ is the only configuration that affects connectivity within an aisle.
For $\Kcal = \{a, b\}$, only the pair $(a,b)$ exists.
This means we only need to consider the connection between $a$ and $b$,
the only way for them to be connected
is by a single edge or through the previous aisle as
seen in Figure \ref{fig:connectivity_sng}.

\input{figures/figure_formulations_single}

\paragraph{Connection (upper bound)}
This constraint says that $a_j$ and $b_j$ can be connected
only if they are connected via a single edge or through the previous aisle.
\begin{equation}
\mbox{\begin{tabular*}{\dimexpr\linewidth-5em\relax}{@{\extracolsep{\fill}}lr@{}}%
$\displaystyle \acr{j}{a}{b} \le \acpass{j}{b} \cond{j>0}{+ \acp{j}{a}{b}}$ &
$\displaystyle j \in \Jcal$%
\end{tabular*}}
\modeltag
\end{equation}

\paragraph{Connected to next aisle (upper bound)}
For each interior aisle $j$, if $a$ has a left edge then
it must be connected to the next aisle,
otherwise the graph is not connected.
This means either it must either directly have
a right edge or be connected to $b$,
which has a right edge
(and symmetrically for $b$).
\begin{equation}
\acr{j}{k}{k'} + \acbarx{j}{k} + \acbarbarx{j}{k}
\ge \acbarx{j-1}{k} + \acbarbarx{j-1}{k} \qquad j \in \Jcal \setminus \{0, m-1\},
(k,k') \in \{ (a,b), (b,a) \}
\modeltag
\end{equation}
\begin{equation}
\acbarx{j}{k'} + \acbarbarx{j}{k'} + \acbarx{j}{k} + \acbarbarx{j}{k}
\ge \acbarx{j-1}{k} + \acbarbarx{j-1}{k}
\qquad j \in \Jcal \setminus \{0, m-1\},
(k,k') \in \{ (a,b), (b,a) \}
\modeltag
\end{equation}

%% file: figures/figure_formulations_single.tex
\begin{figure}[ht]
    \centering
    \begin{subfigure}[t]{0.45\textwidth}
    \centering
    \begin{tikzpicture}[scale=1.2, every node/.style={circle, draw=black, minimum size=6mm}]
      \node (a) at (0, 2) {$a_j$};
      \node (b) at (0, 0)   {$b_j$};
      \draw[thick] (a) -- (-0.9, 2) -- (-0.9, 0) -- (b);
    \end{tikzpicture}
    \caption{Indirect}
    \end{subfigure}
    \hfill
    \begin{subfigure}[t]{0.45\textwidth}
    \centering
    \begin{tikzpicture}[scale=1.2, every node/.style={circle, draw=black, minimum size=6mm}]
      \node (a) at (0, 2) {$a_j$};
      \node (b) at (0, 0)   {$b_j$};
      \draw[thick] (a) -- (b);
    \end{tikzpicture}
    \caption{Direct}
    \end{subfigure}
    \caption[Single-block connectivity paths.]
    {Two ways for crosses $a_j$ and $b_j$ to be connected in aisle $j$.
    (a) Indirectly via the previous aisle.
    (b) Directly though a single traversal.}
    \label{fig:connectivity_sng}
\end{figure}
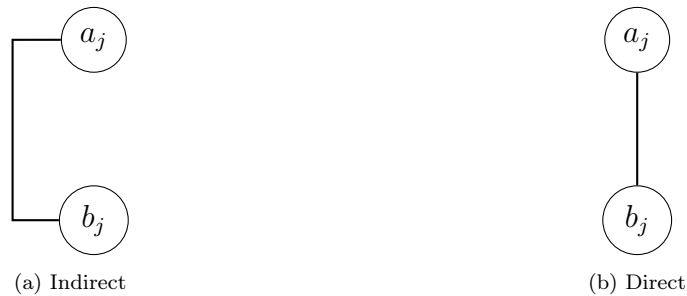

%% file: sections/section_ec_connectivity_twoblock.tex
For a two-block warehouse,
the connection between vertices is more complicated.
Figure \ref{fig:connectivity_two} shows the possible connections for $a_j$, $b_j$, and $c_j$.
The direct and indirect (previous aisle) connections apply again
however, the indorduction of a middle cross-aisle introduces new connections
introduces a number of combinations of connections.

\input{figures/figure_formulations_two}

In this section we will introduce the constraints that enforce these connections.
As the edxtra cross-aisle also complicates the next-aisle connection,
an auxiliary variable is introduced to represent the connection to the next aisle.

\paragraph{Connection of adjacent vertices (upper bound)}
Figure \ref{fig:connectivity_two} shows the possible connections for $a_j$ and $b_j$,
either through the previous aisle, directly (same as single-block), or via $c_j$.
\begin{equation}
\mbox{\begin{tabular*}{\dimexpr\linewidth-5em\relax}{@{\extracolsep{\fill}}lr@{}}%
$\displaystyle \acr{j}{a}{b} \le \acpass{j}{a} + \acp{j}{a}{b} + \acpass{j}{b}$ &
$\displaystyle j \in \Jcal \setminus \{0\}$%
\end{tabular*}}
\modeltag
\end{equation}
\begin{equation}
\mbox{\begin{tabular*}{\dimexpr\linewidth-5em\relax}{@{\extracolsep{\fill}}lr@{}}%
$\displaystyle \acr{j}{a}{b} \le \acpass{j}{a} + \acp{j}{a}{b} + \acp{j}{a}{c}$ &
$\displaystyle j \in \Jcal \setminus \{0\}$%
\end{tabular*}}
\modeltag
\end{equation}

Similarly, for $b_j$ and $c_j$:
\begin{equation}
\mbox{\begin{tabular*}{\dimexpr\linewidth-5em\relax}{@{\extracolsep{\fill}}lr@{}}%
$\displaystyle \acr{j}{b}{c} \le \acpass{j}{b} + \acp{j}{b}{c} + \acpass{j}{c}$ &
$\displaystyle j \in \Jcal \setminus \{0\}$%
\end{tabular*}}
\modeltag
\end{equation}
\begin{equation}
\mbox{\begin{tabular*}{\dimexpr\linewidth-5em\relax}{@{\extracolsep{\fill}}lr@{}}%
$\displaystyle \acr{j}{b}{c} \le \acpass{j}{b} + \acp{j}{b}{c} + \acp{j}{b}{a}$ &
$\displaystyle j \in \Jcal \setminus \{0\}$%
\end{tabular*}}
\modeltag
\end{equation}

\paragraph{Connection of boundary vertices (upper bound)}
Similarly, Figure \ref{fig:connectivity_two} shows the possible connections for $a_j$ and $c_j$,
either through the previous aisle, directly (same as single-block),
or via one of the two combinations including $b_j$.
We simplify this by recognising that the last three
possibilities involve $(a_j, b_j)$ and $(b_j, c_j)$ connections.
\begin{equation}
\mbox{\begin{tabular*}{\dimexpr\linewidth-5em\relax}{@{\extracolsep{\fill}}lr@{}}%
$\displaystyle \acr{j}{a}{c} \le \acr{j}{b}{c} + \acp{j}{a}{c}$ &
$\displaystyle j \in \Jcal \setminus \{0\}$%
\end{tabular*}}
\modeltag
\end{equation}
\begin{equation}
\mbox{\begin{tabular*}{\dimexpr\linewidth-5em\relax}{@{\extracolsep{\fill}}lr@{}}%
$\displaystyle \acr{j}{a}{c} \le \acr{j}{a}{b} + \acp{j}{a}{c}$ &
$\displaystyle j \in \Jcal \setminus \{0\}$%
\end{tabular*}}
\modeltag
\end{equation}

\paragraph{Next-aisle auxiliary (upper bound)}
Each component is must be connected to the next aisle.
We introduce auxiliary variables $\acz{j}{k}{k'}$ to represent
vertex $k$ being connected to the next aisle via vertex $k'$.
This can only happen if $k$ and $k'$ are connected
and $k'$ has a right edge.
\begin{equation}
\mbox{\begin{tabular*}{\dimexpr\linewidth-5em\relax}{@{\extracolsep{\fill}}lr@{}}%
$\displaystyle \acz{j}{k}{k'} \le \acr{j}{k}{k'}$ &
$\displaystyle j \in \Jcal \setminus \{m-1\},\; k, k' \in \Kcal,\; k \neq k'$%
\end{tabular*}}
\modeltag
\end{equation}
\begin{equation}
\mbox{\begin{tabular*}{\dimexpr\linewidth-5em\relax}{@{\extracolsep{\fill}}lr@{}}%
$\displaystyle \acz{j}{k}{k'} \le \acbarx{j}{k'} + \acbarbarx{j}{k'}$ &
$\displaystyle j \in \Jcal \setminus \{m-1\},\; k, k' \in \Kcal,\; k \neq k'$%
\end{tabular*}}
\modeltag
\end{equation}
\begin{equation}
\mbox{\begin{tabular*}{\dimexpr\linewidth-5em\relax}{@{\extracolsep{\fill}}lr@{}}%
$\displaystyle \acz{j}{k}{k'} \in [0,1]$ &
$\displaystyle j \in \Jcal,\; (k,k') \in \mathcal{P}; k < k'$%
\end{tabular*}}
\modeltag
\end{equation}

\paragraph{Next-aisle connectivity}
To enforce this,
if a cross has a left edge,
then it must either have a right edge or be connected to the next aisle
via a different cross.
\begin{equation}
\sum_{k' \in \Kcal: k' \neq k} \acz{j}{k}{k'} + \acbarx{j}{k} + \acbarbarx{j}{k} \ge \acbarx{j-1}{k} + \acbarbarx{j-1}{k} \qquad j \in \Jcal \setminus \{0, m-1\},\; k \in \Kcal
\modeltag
\end{equation}

%% file: figures/figure_formulations_two.tex
\begin{figure}[ht]

    \centering


    \resizebox{\textwidth}{!}{%

    \begin{minipage}{\textwidth}

    \centering

    \begin{tabular}{@{}c@{\hspace{1.25em}}c@{}}

    \begin{tikzpicture}[scale=1, every node/.style={circle, draw=black, minimum size=6mm}]

      \node (a1) at (0, 2) {$a_j$};

      \node (b1) at (0, 0)   {$b_j$};

      \node (c1) at (0, -2) {$c_j$};

      \draw[thick] (a1) -- (-1, 2) -- (-1, 0) -- (b1);

      \node (a2) at (2, 2) {$a_j$};

      \node (b2) at (2, 0)   {$b_j$};

      \node (c2) at (2, -2) {$c_j$};

      \draw[thick] (a2) -- (b2);

      \node (a3) at (4, 2) {$a_j$};

      \node (b3) at (4, 0)   {$b_j$};

      \node (c3) at (4, -2) {$c_j$};

      \draw[thick] (a3) -- (3, 2) -- (3, -2) -- (c3) -- (b3);

    \end{tikzpicture}

    &

    \begin{tikzpicture}[scale=1, every node/.style={circle, draw=black, minimum size=6mm}]

        \node (a1) at (0, 2) {$a_j$};

        \node (b1) at (0, 0)   {$b_j$};

        \node (c1) at (0, -2) {$c_j$};

        \draw[thick] (a1) -- (-1, 2) -- (-1, -2) -- (c1);

        \node (a2) at (2, 2) {$a_j$};

        \node (b2) at (2, 0)   {$b_j$};

        \node (c2) at (2, -2) {$c_j$};

        \draw[thick] (a2) -- (b2) -- (c2);

        \node (a3) at (4, 2) {$a_j$};

        \node (b3) at (4, 0)   {$b_j$};

        \node (c3) at (4, -2) {$c_j$};

        \draw[thick] (a3) -- (3, 2) -- (3, 0) -- (b3) -- (c3);

        \node (a4) at (6, 2) {$a_j$};

        \node (b4) at (6, 0)   {$b_j$};

        \node (c4) at (6, -2) {$c_j$};

        \draw[thick] (a4) -- (b4) -- (5, 0) -- (5, -2) -- (c4);

    \end{tikzpicture}

    \end{tabular}

    \end{minipage}

    }%

    \vspace{0.35\baselineskip}
    \begin{subfigure}[t]{0.47\textwidth}
      \centering
      \subcaption{Connections of $a_j$ and $b_j$.}
    \end{subfigure}\hfill
    \begin{subfigure}[t]{0.47\textwidth}
      \centering
      \subcaption{Connections of $a_j$ and $c_j$.}
    \end{subfigure}

    \medskip
    \caption[Two-block connection paths.]{%
    Adjacent vertices $a_j$ and $b_j$ (left) can be connected as in single-block or via $c_j$
    (symmetrically for $b_j$ and $c_j$).
    Boundary vertices $a_j$ and $c_j$ (right) can be connected through the previous aisle,
    directly, or via one of the two combinations including $b_j$.
    Previous-aisle connections are drawn as a single line for simplicity; they may also involve
    double horizontal edges.%
    }
    \label{fig:connectivity_two}

\end{figure}
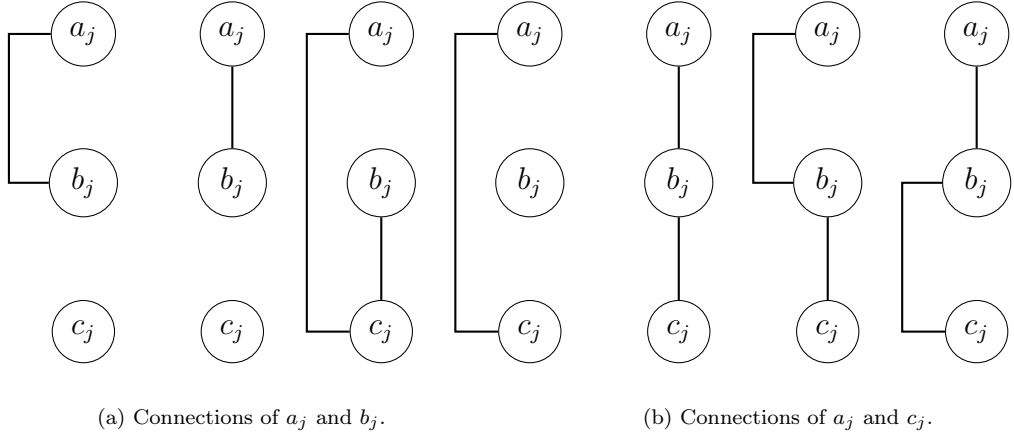

%% file: sections/section_ec_sprpss.tex
For \SPRPsS{}, we use the same notation as \citet{goeke2021modeling}.
Scattered storage augments the EC core with demand coverage and visit-selection constraints; the underlying feasibility semantics (even degree at each vertex, coverage of visited pick positions, and a single connected subgraph containing the depot) are the same as for the standard \SPRP{} connectivity developed above, now applied to the set of positions actually visited.
The algebraic realisation below follows GS and CC but adapts horizontal-configuration and aisle-activation rules to the multi-block setting.

\paragraph{Demand coverage}
This is unchanged from GS and CC,
requiring sufficient supply for each SKU.
\begin{equation}
\mbox{\begin{tabular*}{\dimexpr\linewidth-5em\relax}{@{\extracolsep{\fill}}lr@{}}%
$\displaystyle \sum_{j \in \Jcal} \sum_{i \in \Ical_{jh}} q_{jih} x_{ji} \ge b_h$ &
$\displaystyle h \in \Hcal$%
\end{tabular*}}
\modeltag
\end{equation}

\paragraph{Visit selected positions}
Variable $x_{ji}$ still indicates whether position $(j,i)$ is visited,
however with vertical segments,
we only need to look at the location of interest (and a full traversal),
rather than all positions above and below it in the aisle.
\begin{equation}
\mbox{\begin{tabular*}{\dimexpr\linewidth-5em\relax}{@{\extracolsep{\fill}}lr@{}}%
$\displaystyle \acpass{j}{k} + \pxji{x}{ji} + \qxji{x}{ji} = x_{ji}$ &
$\displaystyle j \in \Jcal,\; i \in \Ical_{jk},\; k \in \Kcal \setminus \{n-1\}$%
\end{tabular*}}
\modeltag
\end{equation}
\begin{equation}
\mbox{\begin{tabular*}{\dimexpr\linewidth-5em\relax}{@{\extracolsep{\fill}}lr@{}}%
$\displaystyle x_{ji} \in \{0,1\}$ &
$\displaystyle j \in \Jcal,\; i \in \Ical_j$%
\end{tabular*}}
\modeltag
\end{equation}

\paragraph{Aisle connectivity}
Define active aisles as $\tilde{x}_j$ as before,
with depot aisle always active.
\begin{equation}
\mbox{\begin{tabular*}{\dimexpr\linewidth-5em\relax}{@{\extracolsep{\fill}}lr@{}}%
$\displaystyle \tilde{x}_j \ge x_{ji}$ &
$\displaystyle j \in \Jcal,\; i \in \Ical_j$%
\end{tabular*}}
\modeltag
\end{equation}
\begin{equation}
\mbox{\begin{tabular*}{\dimexpr\linewidth-5em\relax}{@{\extracolsep{\fill}}l@{}}%
$\displaystyle \tilde{x}_l = 1$%
\end{tabular*}}
\modeltag
\end{equation}
Instead of only considering one horizontal configuration
per aisle,
we now have to look at all cross-aisle configurations.
This requires an upper and lower bound on $\tilde{x}_j$.
\begin{equation}
\mbox{\begin{tabular*}{\dimexpr\linewidth-5em\relax}{@{\extracolsep{\fill}}lr@{}}%
$\displaystyle \acbarx{j}{k} + \acbarbarx{j}{k} \leq \tilde{x}_j$ &
$\displaystyle j \in \Jcal : j < l,\; k \in \Kcal$%
\end{tabular*}}
\modeltag
\end{equation}
\begin{equation}
\mbox{\begin{tabular*}{\dimexpr\linewidth-5em\relax}{@{\extracolsep{\fill}}lr@{}}%
$\displaystyle \sum_{k \in \Kcal} (\acbarx{j}{k} + \acbarbarx{j}{k}) \ge \tilde{x}_j$ &
$\displaystyle j \in \Jcal : j < l$%
\end{tabular*}}
\modeltag
\end{equation}
\begin{equation}
\mbox{\begin{tabular*}{\dimexpr\linewidth-5em\relax}{@{\extracolsep{\fill}}lr@{}}%
$\displaystyle \acbarx{j}{k} + \acbarbarx{j}{k} \leq \tilde{x}_{j+1}$ &
$\displaystyle j \in \Jcal \setminus \{m-1\} : j \ge l,\; k \in \Kcal$%
\end{tabular*}}
\modeltag
\end{equation}
\begin{equation}
\mbox{\begin{tabular*}{\dimexpr\linewidth-5em\relax}{@{\extracolsep{\fill}}lr@{}}%
$\displaystyle \sum_{k \in \Kcal} (\acbarx{j}{k} + \acbarbarx{j}{k}) \ge \tilde{x}_{j+1}$ &
$\displaystyle j \in \Jcal \setminus \{m-1\} : j \ge l$%
\end{tabular*}}
\modeltag
\end{equation}
The propagation of $\tilde{x}_j$ to $\tilde{x}_{j+1}$ is unchanged from GS and CC.
\begin{equation}
\mbox{\begin{tabular*}{\dimexpr\linewidth-5em\relax}{@{\extracolsep{\fill}}lr@{}}%
$\displaystyle \tilde{x}_j \ge \tilde{x}_{j+1}$ &
$\displaystyle j \in \Jcal \setminus \{m-1\} : j \ge l$%
\end{tabular*}}
\modeltag
\end{equation}
\begin{equation}
\mbox{\begin{tabular*}{\dimexpr\linewidth-5em\relax}{@{\extracolsep{\fill}}lr@{}}%
$\displaystyle \tilde{x}_j \le \tilde{x}_{j+1}$ &
$\displaystyle j \in \Jcal : j < l$%
\end{tabular*}}
\modeltag
\end{equation}

\paragraph{Variable domains}
Lastly, we need to adjust
the domain of $\aceta{j}$ to allow for zero edge
pairs in the SPRP-SS formation as not all aisles
must be visited.
\begin{equation}
\mbox{\begin{tabular*}{\dimexpr\linewidth-5em\relax}{@{\extracolsep{\fill}}lr@{}}%
$\displaystyle \aceta{j} \in \{0, \ldots, |\Kcal|\}$ &
$\displaystyle j \in \Jcal \setminus \{m-1\}$%
\end{tabular*}}
\modeltag
\end{equation}

%% file: sections/section_results.tex
\label{sec:computational_experiments}
This section compares the computational performance of the MILP formulations introduced in this paper
on both standard SPRP and SPRP-SS instances.
The goal is to illustrate the practical impact of the structural simplifications and modelling choices.
Since dynamic-programming algorithms provide highly efficient exact methods for the standard SPRP in the
layouts considered here, the main motivation for MILP is to support richer variants such as SPRP-SS
and to enable solver-based integration into larger optimisation settings.
For warehouse analysts and IE practitioners, the experiments therefore emphasize when a structure-aware MILP is preferable to a generic network-flow encoding: smaller presolved models and faster solves translate directly into shorter turnaround for simulation studies or embedded routing calls in higher-level IE optimisation.


\subsection{Experimental setup}

We evaluate the proposed formulations
against existing formulations on randomly generated
SPRP and SPRP-SS instances
for both single-block and two-block layouts.

\paragraph{Benchmarks}
For single-block layouts we compare four formulations:
the edge-based formulation of \citet{goeke2021modeling} (GS),
our Configuration Connectivity formulation (CC),
the network-flow formulation of \citet{hessler2024exact} (NF),
and our Edge Connectivity formulation (EC).
For two-block layouts we compare NF and EC only, as GS and CC are defined for single-block layouts.
NF is included as the state-of-the-art exact method for single and two-block SPRP-SS.

\paragraph{Implementation}
The experiments were implemented in Python on a laptop with
an Intel i7-13800H processor and 32 GB of RAM.
We use Gurobi 11.0.0 \cite{gurobi2024} as the solver.
Unless stated otherwise in the replication materials, we used Gurobi's default parameter settings; every run reported below terminated with optimality proved (zero MIP gap).
The scripts shared for reproduction document any non-default options (for example time limits, thread counts, or presolve flags) if they were changed from the defaults during the study.

\paragraph{Instance generation and parameters}
For standard SPRP we use the same parameter grid:
aisles $m \in \{ 5, 10, 15, 20, 25 \}$
and total pick positions $|P| \in \{ 5, 10, 15, 20, 25 \}$
with 50 instances per combination,
giving 1250 instances in total (no scattered-storage parameter).
SPRP-SS instances are generated following \citet{goeke2021modeling}.
We vary the number of aisles and the number of requested SKUs (articles),
and we consider multiple
values of the scattered-storage parameter $\alpha$.
Benchmark instances were generated for all combinations of
$\alpha \in \{1, 2, 3, 4, 5 \}$ (only for SPRP-SS),
aisles $m \in \{ 5, 10, 15, 20, 25 \}$
and total pick positions
$|P| \in \{ 5, 10, 15, 20, 25 \}$
with 50 instances per combination,
giving 6250 SPRP-SS instances.
The depot was randomly assigned to an aisle,
equally likely in the front or rear cross-aisle.
The number of pick locations per aisle was fixed at 90.
We chose not to vary the number of pick locations per aisle
as all formulations are independent of the number of pick locations.

We follow the instance generation procedure of \citet{goeke2021modeling}:
the number of distinct SKUs in the warehouse is set
to $\xi = \max(a, \lceil m \cdot n / \alpha \rceil)$,
where $a$ is the number of SKUs in the pick list,
$m \cdot n$ is the storage capacity,
and the scatter factor $\alpha \ge 1$ controls the degree of duplication.
For $\alpha = 1$, each picking position holds a different SKU (no duplication).
A larger $\alpha$ yields fewer distinct SKUs and more positions per SKU,
so the solver has more flexibility in where to pick each requested item
and the problem becomes
more representative of scattered-storage settings.
SKUs and pick lists are assigned to storage positions
and demand according to turnover-based classes
as in \citet{goeke2021modeling},
we use our own parameter grid ($\alpha$, $m$, $|P|$) as above.
The instance-generation code supplied with the replication package records pseudorandom draws (including any fixed seeds) so that the same benchmark sets can be regenerated.


\subsection{Single-Block Results}
\label{ssec:singleblock_results}

We compare four formulations (GS, CC, NF, EC) for single-block layouts.
Results are given for standard SPRP first, then for SPRP-SS.

\input{tables/table_singleblock_SPRP_overall}
Table~\ref{tab:singleblock_SPRP_overall} presents the overall performance metrics
for single-block formulations (GS, CC, NF, EC) on standard SPRP instances.
All four formulations solve all 1250 instances to optimality.
Although standard SPRP is not the main focus of this paper
as more efficient dynamic programming algorithms exist,
the table illustrates the relative runtimes aligning with expectations.
NF is by far the fastest (on all overall metrics),
as expected since it encodes the efficient dynamic-programming structure directly
as a shortest-path network flow.
CC shows improved runtimes over GS,
indicating that the structural simplifications and reduced variable count are effective.
EC shows slightly worse runtimes than GS and CC,
aligning with the additional variables and constraints required
to enforce connectivity in a general setting.
The differences are small, however,
indicating that the additional complexity is not significant.

\input{tables/table_singleblock_overall}

Table~\ref{tab:singleblock_overall} summarizes overall performance for single-block SPRP-SS,
and Figure~\ref{plt:singleblock_combined} shows runtimes as a function of $\alpha$ (top panel),
number of aisles (middle panel), and number of articles (bottom panel), all on a logarithmic scale.
All formulations solve all 6250 SPRP-SS instances to optimality.
On average, CC is the best performer,
consistently best on all metrics.
EC shows similar behaviors to CC over the parameter space,
with slightly worse runtimes.
The same pattern is observed for GS with consistently slower runtimes across the parameter space.
Over the $\alpha$ parameter, all four formulations scale similarly, including NF.
On the other two parameters, the differences of NF are more pronounced.
The number of aisles has a significant impact on the three arc-based formulations,
however, NF remains roughly flat.
On the other hand,
in contrast to the number of articles not impacting the other three,
it has a significant impact on NF,
with significantly worse runtimes as the number of articles increases.
Observing Table~\ref{tab:singleblock_overall},
we see that although NF has the worst average runtime
and similar median to GS,
it has the second best geometric mean,
indicating that the particular sensitivity to the
number of articles has skewed the other metrics
towards worse values.

\input{plots/plot_singleblock_combined}

For two-block layouts we compare only NF and EC below, as GS and CC are defined for single-block layouts only.


\subsection{Two-Block Results}
\label{ssec:twoblock_results}

For two-block layouts we compare NF and EC on the same instance sets and metrics as above.

\input{tables/table_twoblock_SPRP_overall}

Table~\ref{tab:twoblock_SPRP_overall} presents the overall performance metrics
for two-block formulations (NF, EC) on standard SPRP instances.
All 1250 instances are solved to optimality.
NF is much faster than EC on all metrics,
again consistent with the network-flow formulation matching the DP structure.

\input{tables/table_twoblock_overall}

Table~\ref{tab:twoblock_overall} gives overall metrics for two-block SPRP-SS,
and Figure~\ref{plt:twoblock_combined} compares runtimes for NF and EC with varying $\alpha$ (top),
number of aisles (middle), and number of articles (bottom), on a logarithmic scale.
All 6250 instances are solved to optimality.
On both standard and geometric means, NF is faster than EC,
however EC has a better median runtime,
indicating a sensitivity to the more difficult instances.
By parameter, both show similar behavior over the $\alpha$ parameter,
with NF consistently faster.
Over the number of aisles,
NF exhibits a flat behavior,
while EC shows a significant increase in runtime as the number of aisles increases,
performing much better than NF on instances with a small number of aisles
and much worse on instances with a large number of aisles.
Contrastingly, EC shows a flat behavior over the number of articles,
with NF performing much better on instances with a small number of articles
and degrading as the number of articles increases,
performing worse than EC on instances with a large number of articles.

\input{plots/plot_twoblock_combined}


\subsection{Analysis of Results}
\label{ssec:analysis_results}

We synthesise the single-block and two-block results and interpret them in light
of the formulations' structure.

\paragraph{Standard SPRP (single-block and two-block)}
As shown in Table~\ref{tab:singleblock_SPRP_overall} and Table~\ref{tab:twoblock_SPRP_overall},
for standard SPRP both single-block and two-block formulations solve all instances to optimality.
NF's advantage is explained by its direct encoding of the DP state graph
as a shortest-path network flow, which is an efficient way to ensure connectivity.
GS, CC, and EC are graph-based formulations that exploit the DP structure to some degree
but do not encode it as compactly.
The NF method is therefore the more efficient choice for standard SPRP
when a MILP is used.

\paragraph{Single-block SPRP-SS}
The single-block SPRP-SS results (Table~\ref{tab:singleblock_overall} and Figure~\ref{plt:singleblock_combined})
show that CC consistently outperforms GS, as expected:
CC implements structural simplifications
(no double traversals, connectivity in terms of horizontal configurations only)
and thus has fewer variables and simpler constraints.
EC shows slightly worse runtimes than CC, also as expected,
because it requires additional variables and constraints to enforce connectivity in a general setting
that supports both single-block and two-block layouts.
EC is also slightly worse than GS on standard SPRP for the same reasons;
on SPRP-SS, however, EC is slightly better than GS,
which may be explained by EC modelling the warehouse structure more completely
through the relationship between adjacent vertices and tighter bounding constraints.
GS, CC, and EC show similar behaviour over all parameters, as is expected for arc-based formulations.
Over the $\alpha$ parameter, all four formulations scale similarly, including NF.
NF's sensitivity to the number of articles is consistent with a formulation
whose size grows with the number of items;
as the number of aisles increases, NF remains roughly flat while the three graph-based formulations
show a significant increase in runtime.

\paragraph{Two-block SPRP-SS}
The two-block SPRP-SS results (Table~\ref{tab:twoblock_overall} and Figure~\ref{plt:twoblock_combined})
reflect a trade-off between the two formulations.
NF performs better overall and scales well with the number of aisles,
as it encodes the DP state graph directly as a shortest-path network flow
and considers horizontal edge configurations as a whole.
EC reasons at the level of individual vertices to generalise between single-block and two-block layouts
(and potentially more blocks in the future),
which adds flexibility but cost.
NF is sensitive to the number of articles, however,
and EC performs better than NF on instances with a large number of articles,
as well as better on instances with a small number of aisles.
That EC can outperform NF in these regimes shows that the graph-based modelling of the warehouse structure
and the connection constraints of the EC formulation successfully exploit the layout.
The trade-off is practical: EC is preferable when the number of aisles is small
or the number of requested articles is large,
whereas NF is preferable when the number of requested articles is small;
as parameters grow, the choice depends on which dimension dominates.

\paragraph{Summary}
The results are consistent with the formulation design goals:
simplifying the edge-based model by removing redundant traversal/configuration options and streamlining connectivity
(CC) yields clear runtime gains over the baseline (GS),
while the generalised formulation (EC) pays a moderate overhead for layout independence.
For practitioners, the choice of formulation depends on the layout (single-block vs.\ two-block)
and on instance dimensions (number of aisles and number of articles),
with NF strongest for standard SPRP and for SPRP-SS when the number of articles is small,
and EC competitive or preferable when the number of aisles is small or the number of articles is large.
This is consistent with our understanding of the formulations as NF models
the state transitions of a valid tour efficiently but grows quadratically
with the number of items in an aisle due to the combinations
of gap configurations considered \cite{luke2025linear}.
On the other hand,
the edge-based models are linear in the number of items but
are restricted by the complexity of modelling connectivity.

%% file: tables/table_singleblock_SPRP_overall.tex
\begin{table}[t]
\centering
\caption{
    Overall performance metrics for single-block warehouse formulations
    (GS, CC, NF, EC), standard SPRP, over all test instances.
}
\label{tab:singleblock_SPRP_overall}
\vspace{1em}
\begin{filecontents*}{plots/singleblock_SPRP_results_overall.csv}
metric,GS,CC,NF,EC
solved (1250),1250,1250,1250,1250
overall average,7.40,4.96,1.69,8.00
median,5.99,4.30,1.60,6.37
geometric mean,6.26,4.35,1.52,6.57
\end{filecontents*}
\pgfplotstabletypeset[
    col sep=comma,
    string type,
    columns={metric,GS,CC,NF,EC},
    columns/metric/.style={column name=Metric, column type=l},
    columns/GS/.style={column name=GS, column type=c},
    columns/CC/.style={column name=CC, column type=c},
    columns/NF/.style={column name=NF, column type=c},
    columns/EC/.style={column name=EC, column type=c},
    every head row/.style={
        before row=\hline,
        after row=\hline
    },
    every last row/.style={after row=\hline}
]{plots/singleblock_SPRP_results_overall.csv}
\end{table}

%% file: tables/table_singleblock_overall.tex
\begin{table}[t]
\centering
\caption{
    Overall performance metrics for single-block warehouse formulations
    (GS, NF, CC, EC), SPRP-SS, over all test instances.
}
\label{tab:singleblock_overall}
\vspace{1em}
\begin{filecontents*}{plots/singleblock_SPRP-SS_results_overall.csv}
metric,GS,CC,NF,EC
solved (6250),6250,6250,6250,6250
overall average,116.79,76.97,124.57,81.12
median,49.85,34.93,48.48,40.51
geometric mean,47.02,32.69,38.80,40.27
\end{filecontents*}
\pgfplotstabletypeset[
    col sep=comma,
    string type,
    columns={metric,GS,CC,NF,EC},
    columns/metric/.style={column name=Metric, column type=l},
    columns/GS/.style={column name=GS, column type=c},
    columns/CC/.style={column name=CC, column type=c},
    columns/NF/.style={column name=NF, column type=c},
    columns/EC/.style={column name=EC, column type=c},
    every head row/.style={
        before row=\hline,
        after row=\hline
    },
    every last row/.style={after row=\hline}
]{plots/singleblock_SPRP-SS_results_overall.csv}
\end{table}

%% file: plots/plot_singleblock_combined.tex
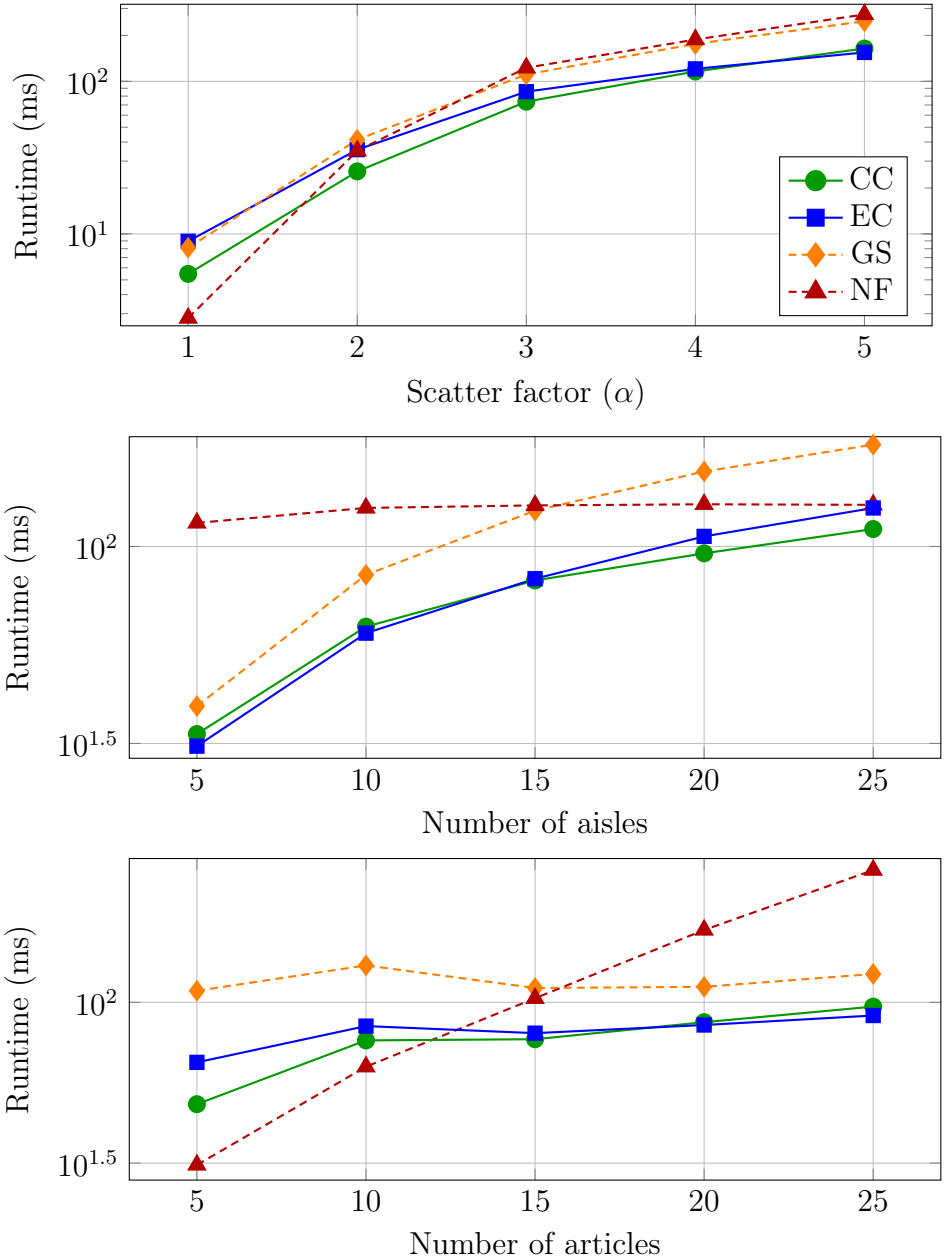
\begin{figure}[t]
    \centering
    \begin{subfigure}[b]{0.95\textwidth}
        \centering
        \begin{tikzpicture}
        \begin{axis}[
            width=0.95\linewidth,
            height=0.45\linewidth,
            xlabel={Scatter factor ($\alpha$)},
            ylabel={Runtime (\si{\milli\second})},
            legend pos=south east,
            legend image post style={scale=1.2},
            legend style={font=\normalsize},
            grid=major,
            xtick={1,2,3,4,5},
            ymode=log,
            ymin=2.5,
            ymax=320,
            mark size=3pt
        ]

        \addplot[color=green!60!black, thick, mark=*, mark options={solid}] table[col sep=comma, x=alpha, y=CC]{plots/singleblock_SPRP-SS_results_alpha.csv};
        \addlegendentry{CC}
        
        \addplot[blue, thick, mark=square*, mark size=2.5pt, mark options={solid}] table[col sep=comma, x=alpha, y=EC]{plots/singleblock_SPRP-SS_results_alpha.csv};
        \addlegendentry{EC}
        
        \addplot[color=orange, thick, dash pattern=on 3pt off 2pt, mark=diamond*, mark size=3.5pt, mark options={solid}] table[col sep=comma, x=alpha, y=GS]{plots/singleblock_SPRP-SS_results_alpha.csv};
        \addlegendentry{GS}
        
        \addplot[color=red!70!black, thick, densely dashed, mark=triangle*, mark size=3.5pt, mark options={solid}] table[col sep=comma, x=alpha, y=NF]{plots/singleblock_SPRP-SS_results_alpha.csv};
        \addlegendentry{NF}
        
        \end{axis}
        \end{tikzpicture}
    \end{subfigure}
    
    \vspace{0.4em}
    
    \begin{subfigure}[b]{0.95\textwidth}
        \centering
        \begin{tikzpicture}
        \begin{axis}[
            width=0.95\linewidth,
            height=0.45\linewidth,
            xlabel={Number of aisles},
            ylabel={Runtime (\si{\milli\second})},
            grid=major,
            xtick={5,10,15,20,25},
            ymode=log,
            ymin=29,
            ymax=190,
            mark size=3pt
        ]
        
        \addplot[color=orange, thick, dash pattern=on 3pt off 2pt, mark=diamond*, mark size=3.5pt, mark options={solid}] table[col sep=comma, x=num_aisles, y=GS]{plots/singleblock_SPRP-SS_results_aisles.csv};
        
        \addplot[color=red!70!black, thick, densely dashed, mark=triangle*, mark size=3.5pt, mark options={solid}] table[col sep=comma, x=num_aisles, y=NF]{plots/singleblock_SPRP-SS_results_aisles.csv};
        
        \addplot[color=green!60!black, thick, mark=*, mark options={solid}] table[col sep=comma, x=num_aisles, y=CC]{plots/singleblock_SPRP-SS_results_aisles.csv};

        \addplot[blue, thick, mark=square*, mark size=2.5pt, mark options={solid}] table[col sep=comma, x=num_aisles, y=EC]{plots/singleblock_SPRP-SS_results_aisles.csv};
        
        \end{axis}
        \end{tikzpicture}
    \end{subfigure}
    
    \vspace{0.4em}
    
    \begin{subfigure}[b]{0.95\textwidth}
        \centering
        \begin{tikzpicture}
        \begin{axis}[
            width=0.95\linewidth,
            height=0.45\linewidth,
            xlabel={Number of articles},
            ylabel={Runtime (\si{\milli\second})},
            grid=major,
            xtick={5,10,15,20,25},
            ymode=log,
            ymin=28,
            ymax=280,
            mark size=3pt
        ]
        
        \addplot[color=orange, thick, dash pattern=on 3pt off 2pt, mark=diamond*, mark size=3.5pt, mark options={solid}] table[col sep=comma, x=num_articles, y=GS]{plots/singleblock_SPRP-SS_results_articles.csv};
        
        \addplot[color=red!70!black, thick, densely dashed, mark=triangle*, mark size=3.5pt, mark options={solid}] table[col sep=comma, x=num_articles, y=NF]{plots/singleblock_SPRP-SS_results_articles.csv};
        
        \addplot[color=green!60!black, thick, mark=*, mark options={solid}] table[col sep=comma, x=num_articles, y=CC]{plots/singleblock_SPRP-SS_results_articles.csv};

        \addplot[blue, thick, mark=square*, mark size=2.5pt, mark options={solid}] table[col sep=comma, x=num_articles, y=EC]{plots/singleblock_SPRP-SS_results_articles.csv};
        
        \end{axis}
        \end{tikzpicture}
    \end{subfigure}
    
\caption[Runtime comparison of single-block formulations against baselines.]
{
    Comparison of formulation runtimes (log scale) for single-block warehouses
    with varying parameters: scatter factor ($\alpha$) (top), number of aisles (middle), and number of articles (bottom).
}
\label{plt:singleblock_combined}
\end{figure}

%% file: tables/table_twoblock_SPRP_overall.tex
\begin{table}[t]
\centering
\caption{
    Overall performance metrics for two-block warehouse formulations
    (NF, EC), standard SPRP, over all test instances.
}
\label{tab:twoblock_SPRP_overall}
\vspace{1em}
\begin{filecontents*}{plots/twoblock_SPRP_results_overall.csv}
metric,NF,EC
solved (1250),1250,1250
overall average,14.40,106.02
median,12.84,72.64
geometric mean,12.13,62.78
\end{filecontents*}
\pgfplotstabletypeset[
    col sep=comma,
    string type,
    columns={metric,NF,EC},
    columns/metric/.style={column name=Metric, column type=l},
    columns/NF/.style={column name=NF, column type=c},
    columns/EC/.style={column name=EC, column type=c},
    every head row/.style={
        before row=\hline,
        after row=\hline
    },
    every last row/.style={after row=\hline}
]{plots/twoblock_SPRP_results_overall.csv}
\end{table}

%% file: tables/table_twoblock_overall.tex
\begin{table}[t]
\centering
\caption{
    Overall performance metrics for two-block warehouse formulations
    (NF, EC), SPRP-SS, over all test instances.
}
\label{tab:twoblock_overall}
\vspace{1em}
\begin{filecontents*}{plots/twoblock_SPRP-SS_results_overall.csv}
metric,NF,EC
solved (6250),6250,6250
overall average,1297.64,1842.06
median,424.24,352.83
geometric mean,344.57,397.12
\end{filecontents*}
\pgfplotstabletypeset[
    col sep=comma,
    string type,
    columns={metric,NF,EC},
    columns/metric/.style={column name=Metric, column type=l},
    columns/NF/.style={column name=NF, column type=c},
    columns/EC/.style={column name=EC, column type=c},
    every head row/.style={
        before row=\hline,
        after row=\hline
    },
    every last row/.style={after row=\hline}
]{plots/twoblock_SPRP-SS_results_overall.csv}
\end{table}

%% file: plots/plot_twoblock_combined.tex
\begin{figure}[t]
    \centering
    \begin{subfigure}[b]{0.95\textwidth}
        \centering
        \begin{tikzpicture}
        \begin{axis}[
            width=0.95\linewidth,
            height=0.45\linewidth,
            xlabel={Scatter factor ($\alpha$)},
            ylabel={Runtime (\si{\milli\second})},
            legend pos=south east,
            legend image post style={scale=1.2},
            legend style={font=\normalsize},
            grid=major,
            xtick={1,2,3,4,5},
            ymode=log,
            ymin=15,
            ymax=4800,
            mark size=3pt
        ]

        \addplot[blue, thick, mark=square*, mark size=2.5pt, mark options={solid}] table[col sep=comma, x=alpha, y=EC]{plots/twoblock_SPRP-SS_results_alpha.csv};
        \addlegendentry{EC}
        
        \addplot[color=red!70!black, thick, densely dashed, mark=triangle*, mark size=3.5pt, mark options={solid}] table[col sep=comma, x=alpha, y=NF]{plots/twoblock_SPRP-SS_results_alpha.csv};
        \addlegendentry{NF}
        
        \end{axis}
        \end{tikzpicture}
    \end{subfigure}
    
    \vspace{0.4em}
    
    \begin{subfigure}[b]{0.95\textwidth}
        \centering
        \begin{tikzpicture}
        \begin{axis}[
            width=0.95\linewidth,
            height=0.45\linewidth,
            xlabel={Number of aisles},
            ylabel={Runtime (\si{\milli\second})},
            grid=major,
            xtick={5,10,15,20,25},
            ymode=log,
            mark size=3pt
        ]

        \addplot[blue, thick, mark=square*, mark size=2.5pt, mark options={solid}] table[col sep=comma, x=num_aisles, y=EC]{plots/twoblock_SPRP-SS_results_aisles.csv};
        
        \addplot[color=red!70!black, thick, densely dashed, mark=triangle*, mark size=3.5pt, mark options={solid}] table[col sep=comma, x=num_aisles, y=NF]{plots/twoblock_SPRP-SS_results_aisles.csv};
                
        \end{axis}
        \end{tikzpicture}
    \end{subfigure}
    
    \vspace{0.4em}
    
    \begin{subfigure}[b]{0.95\textwidth}
        \centering
        \begin{tikzpicture}
        \begin{axis}[
            width=0.95\linewidth,
            height=0.45\linewidth,
            xlabel={Number of articles},
            ylabel={Runtime (\si{\milli\second})},
            grid=major,
            xtick={5,10,15,20,25},
            ymode=log,
            mark size=3pt
        ]
        
        \addplot[color=red!70!black, thick, densely dashed, mark=triangle*, mark size=3.5pt, mark options={solid}] table[col sep=comma, x=num_articles, y=NF]{plots/twoblock_SPRP-SS_results_articles.csv};
        
        \addplot[blue, thick, mark=square*, mark size=2.5pt] table[col sep=comma, x=num_articles, y=EC, mark options={solid}]{plots/twoblock_SPRP-SS_results_articles.csv};
        
        \end{axis}
        \end{tikzpicture}
    \end{subfigure}
    
\caption[Runtime comparison of two-block formulation against baseline.]
{
    Comparison of formulation runtimes (log scale) for two-block warehouses
    with varying parameters: scatter factor ($\alpha$) (top), number of aisles (middle), and number of articles (bottom).
}
\label{plt:twoblock_combined}
\end{figure}
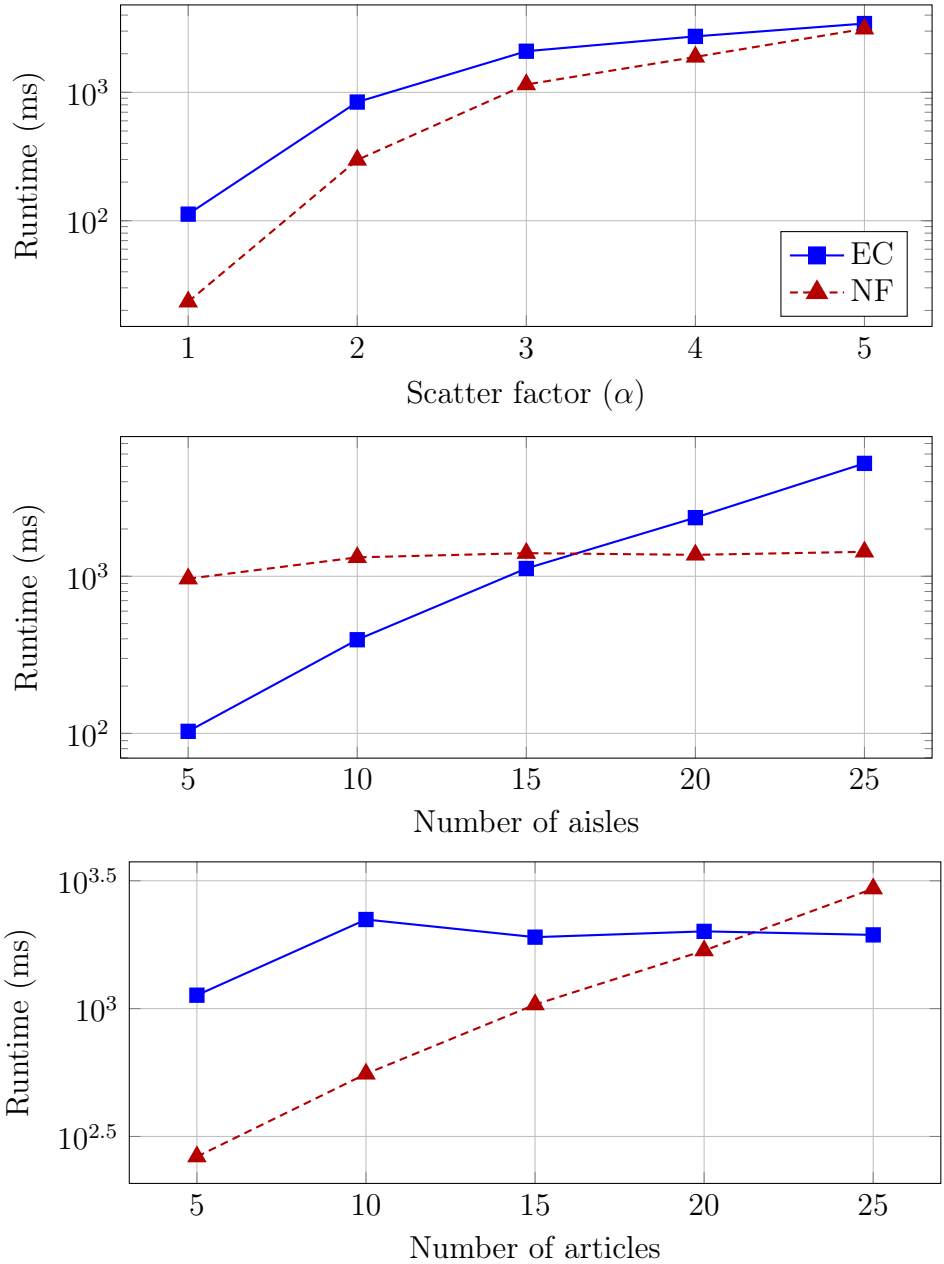

%% file: sections_conclusion.tex
We developed two structure-aware MILPs for the \SPRP{} and \SPRPsS{} in rectangular warehouses: \emph{Configuration Connectivity} (CC), a single-block formulation that strengthens the edge-based baseline of \citet{goeke2021modeling} (GS) by removing redundant configurations and simplifying connectivity, and \emph{Edge Connectivity} (EC), a graph-based model with a common core and layout-specific connection rules for single- and two-block designs.
Both extend to scattered storage.
Computational experiments on large random instance sets compared CC and EC to GS and to the network-flow formulation of \citet{hessler2024exact} (NF), quantifying when the structural restrictions translate into faster solves or smaller presolved models.

\subsection{Implications for practice}

Industrial engineers and OR analysts who embed routing inside warehouse planning tools face a choice between dedicated dynamic-programming codes, network-flow encodings of the DP state graph (NF), and arc-based MILPs (GS, CC, EC).
For \textbf{single-block} layouts, CC is a practical default when a \textbf{compact} MILP is needed for repeated calls (e.g.\ inside slotting, storage-policy, or layout what-if models).
For \textbf{two-block} layouts, EC offers one implementation path for both block types; relative to NF it can be preferable when the number of aisles is small or the number of line items is large, whereas NF is often stronger when pick lists are short, matching the structural trade-off between state-graph compactness and article-driven growth documented in Section~\ref{sec:computational_experiments}.
Where \textbf{scattered storage} applies or routing must sit alongside \textbf{batching, workforce, or design} constraints, solver-based MILPs remain attractive because they avoid maintaining bespoke DP machinery while still delivering proven optimal tours at the scales we tested.

\subsection{Limitations}

The theory and experiments concern \textbf{rectangular, block-structured} layouts and \textbf{rectilinear} travel; fishbone, chevron, or free-form floors are not covered.
We treat \textbf{static, single-picker} routing from a \textbf{single depot}; dynamic order arrival, online replanning, and coordinated multi-picker systems are left open.
Although EC is designed to generalise beyond two blocks, we only report computational evidence for \textbf{single- and two-block} instances.

\subsection{Future research}

Extending EC connectivity to \textbf{three or more cross-aisles} would widen industrial applicability and requires tractable connection descriptions across additional block interfaces.
\textbf{Operational extensions}, such as multiple depots, time windows, congestion, or handling rules, are natural candidates for the same MILP shell because they attach cleanly to a graph-based routing core.
\textbf{Integration} with batching, sequencing, and storage assignment remains a priority for IE research: the formulations here can act as exact routing subroutines or cost oracles inside larger decision models once those upper-level problems are specified.

%% file: sections/appendix_gs.tex
\section{Goeke--Schneider formulation (GS)}
\label{sec:appendix_gs}

This appendix records the complete GS formulation (as used as a baseline for CC in the thesis sources),
kept complete for reference.
Only minimal wording has been adjusted to remove thesis-only cross references.

For an illustration of the vertical configuration variables used in GS, see Figure~\ref{fig:formulations_gs}.

\subsection{Standard SPRP formulation}

\setcounter{equation}{0}
\paragraph{Objective (Minimize tour length)}
\begin{equation}
\min \sum_{j \in \Jcal} \xx{c}_j \xx{x}_j + \xb{c}_j \xb{x}_j + \xbx{c}_j \xbx{x}_j + \xxb{c}_j \xxb{x}_j + c^{\vert}_j \xpass{j} + c^{\vert\vert}_j \xppass{j} + \sum_{j \in \Jcal} \sum_{i \in \Ical_j} \left( \pxji{c}{ji} \pxji{x}{ji} + \qxji{c}{ji} \qxji{x}{ji} \right)
\end{equation}

\paragraph{Cross-aisle configurations}
\begin{equation}
\mbox{\begin{tabular*}{\dimexpr\linewidth-5em\relax}{@{\extracolsep{\fill}}lr@{}}
$\displaystyle \xx{x}_j + \xb{x}_j + \xbx{x}_j + \xxb{x}_j = 1$ &
$\displaystyle j \in \Jcal \setminus \{m-1\}$
\end{tabular*}}
\end{equation}

\paragraph{Visit all required positions}
\begin{equation}
\mbox{\begin{tabular*}{\dimexpr\linewidth-5em\relax}{@{\extracolsep{\fill}}lr@{}}
$\displaystyle \xpass{j} + \xppass{j} + \sum_{i' \in \Ical_j : i' \geq i} \qxji{x}{ji'} + \sum_{i' \in \Ical_j : i' \leq i} \pxji{x}{ji'} \geq 1$ &
$\displaystyle j \in \Jcal,\; i \in \Ical_j$
\end{tabular*}}
\end{equation}

\paragraph{Top and bottom connectivity}
\begin{equation}
\cond{j>0}{\xx{x}_{j-1} + \xbx{x}_{j-1} + \xxb{x}_{j-1}} + \xx{x}_j + \xbx{x}_j + \xxb{x}_j \geq \pxji{x}{ji} \qquad \begin{cases} \text{if } \thetad = 1: \; j \in \Jcal \\ \text{else: } j \in \Jcal \setminus \{l\} \end{cases},\; i \in \Ical_j
\end{equation}
\begin{equation}
\cond{j>0}{\xb{x}_{j-1} + \xbx{x}_{j-1} + \xxb{x}_{j-1}} + \xb{x}_j + \xbx{x}_j + \xxb{x}_j \geq \qxji{x}{ji} \qquad \begin{cases} \text{if } \thetad = 0: \; j \in \Jcal \\ \text{else: } j \in \Jcal \setminus \{l\} \end{cases},\; i \in \Ical_j
\end{equation}

\paragraph{Switches between top and bottom}
\begin{equation}
\mbox{\begin{tabular*}{\dimexpr\linewidth-5em\relax}{@{\extracolsep{\fill}}lr@{}}
$\displaystyle \xx{x}_{j-1} + \xb{x}_j \leq \xppass{j} + 1$ &
$\displaystyle j \in \Jcal \setminus \{0\}$
\end{tabular*}}
\end{equation}
\begin{equation}
\mbox{\begin{tabular*}{\dimexpr\linewidth-5em\relax}{@{\extracolsep{\fill}}lr@{}}
$\displaystyle \xb{x}_{j-1} + \xx{x}_j \leq \xppass{j} + 1$ &
$\displaystyle j \in \Jcal \setminus \{0\}$
\end{tabular*}}
\end{equation}

\paragraph{Depot inclusion}
\begin{equation}
\mbox{\begin{tabular*}{\dimexpr\linewidth-5em\relax}{@{\extracolsep{\fill}}lr@{}}
$\displaystyle 2\xppass{l} + \xpass{l} + \cond{l>0}{\xx{x}_{l-1} + \xxb{x}_{l-1}} + \xb{x}_l + \xxb{x}_l \geq \cond{l>0}{\xb{x}_{l-1}} + \xx{x}_l$ &
$\displaystyle \text{if } \thetad = 1$
\end{tabular*}}
\end{equation}
\begin{equation}
\mbox{\begin{tabular*}{\dimexpr\linewidth-5em\relax}{@{\extracolsep{\fill}}lr@{}}
$\displaystyle 2\xppass{l} + \xpass{l} + \cond{l>0}{\xb{x}_{l-1} + \xxb{x}_{l-1}} + \xx{x}_l + \xxb{x}_l \geq \cond{l>0}{\xx{x}_{l-1}} + \xb{x}_l$ &
$\displaystyle \text{if } \thetad = 0$
\end{tabular*}}
\end{equation}

\paragraph{Degree parity}
\begin{equation}
\mbox{\begin{tabular*}{\dimexpr\linewidth-5em\relax}{@{\extracolsep{\fill}}lr@{}}
$\displaystyle \cond{j>0}{\xbx{x}_{j-1} + 2\xxb{x}_{j-1} + 2\xb{x}_{j-1}} + \xbx{x}_j + 2\xxb{x}_j + 2\xb{x}_j + \xpass{j} + 2\xppass{j} = 2\pitop{j}$ &
$\displaystyle j \in \Jcal$
\end{tabular*}}
\end{equation}
\begin{equation}
\mbox{\begin{tabular*}{\dimexpr\linewidth-5em\relax}{@{\extracolsep{\fill}}lr@{}}
$\displaystyle \cond{j>0}{\xbx{x}_{j-1} + 2\xxb{x}_{j-1} + 2\xx{x}_{j-1}} + \xbx{x}_j + 2\xxb{x}_j + 2\xx{x}_j + \xpass{j} + 2\xppass{j} = 2\pibot{j}$ &
$\displaystyle j \in \Jcal$
\end{tabular*}}
\end{equation}

\paragraph{Connected components}
\begin{equation}
\mbox{\begin{tabular*}{\dimexpr\linewidth-5em\relax}{@{\extracolsep{\fill}}lr@{}}
$\displaystyle \xxb{x}_j + \xx{x}_{j-1} + \xb{x}_{j-1} - \xppass{j} \leq \tau_j + 1$ &
$\displaystyle j \in \Jcal \setminus \{0\}$
\end{tabular*}}
\end{equation}

\paragraph{Component propagation when left part unvisited}
\begin{equation}
\mbox{\begin{tabular*}{\dimexpr\linewidth-5em\relax}{@{\extracolsep{\fill}}lr@{}}
$\displaystyle \xxb{x}_j + \cond{j>0}{-\xxb{x}_{j-1} - \xx{x}_{j-1} - \xb{x}_{j-1}} - \xppass{j} - \xpass{j} \leq \tau_j$ &
$\displaystyle j \in \Jcal$
\end{tabular*}}
\end{equation}

\paragraph{Propagate components}
\begin{equation}
\mbox{\begin{tabular*}{\dimexpr\linewidth-5em\relax}{@{\extracolsep{\fill}}lr@{}}
$\displaystyle \tau_{j-1} - \xpass{j} - \xppass{j} \leq \tau_j$ &
$\displaystyle j \in \Jcal \setminus \{0\}$
\end{tabular*}}
\end{equation}

\paragraph{Ensure 22 when two components}
\begin{equation}
\mbox{\begin{tabular*}{\dimexpr\linewidth-5em\relax}{@{\extracolsep{\fill}}lr@{}}
$\displaystyle \tau_j \leq \xxb{x}_j$ &
$\displaystyle j \in \Jcal$
\end{tabular*}}
\end{equation}

\paragraph{Variable domains}
\begin{equation}
\mbox{\begin{tabular*}{\dimexpr\linewidth-5em\relax}{@{\extracolsep{\fill}}lr@{}}
$\displaystyle \xx{x}_j, \xb{x}_j, \xbx{x}_j, \xxb{x}_j, \tau_j \in \{0,1\}$ &
$\displaystyle j \in \Jcal \setminus \{m-1\}$
\end{tabular*}}
\end{equation}
\begin{equation}
\mbox{\begin{tabular*}{\dimexpr\linewidth-5em\relax}{@{\extracolsep{\fill}}lr@{}}
$\displaystyle \xpass{j}, \xppass{j} \in \{0,1\}$ &
$\displaystyle j \in \Jcal$
\end{tabular*}}
\end{equation}
\begin{equation}
\mbox{\begin{tabular*}{\dimexpr\linewidth-5em\relax}{@{\extracolsep{\fill}}lr@{}}
$\displaystyle \pxji{x}{ji}, \qxji{x}{ji} \in \{0,1\}$ &
$\displaystyle j \in \Jcal,\; i \in \Ical_j$
\end{tabular*}}
\end{equation}
\begin{equation}
\mbox{\begin{tabular*}{\dimexpr\linewidth-5em\relax}{@{\extracolsep{\fill}}lr@{}}
$\displaystyle \pitop{j}, \pibot{j} \in \Ncal$ &
$\displaystyle j \in \Jcal$
\end{tabular*}}
\end{equation}

\paragraph{Last aisle}
\begin{flalign}
& \xx{x}_{m-1} = \xb{x}_{m-1} = \xbx{x}_{m-1} = \xxb{x}_{m-1} = \tau_{m-1} = 0 &
\end{flalign}

\subsection{Adapting to SPRP with scattered storage}

Replace (3) by (21)--(23). Replace (2) by (24)--(30).

\paragraph{Demand coverage}
\begin{equation}
\mbox{\begin{tabular*}{\dimexpr\linewidth-5em\relax}{@{\extracolsep{\fill}}lr@{}}%
$\displaystyle \sum_{j \in \Jcal} \sum_{i \in \Ical_{jh}} q_{jih} x_{ji} \geq b_h$ &
$\displaystyle h \in \Hcal$%
\end{tabular*}}
\end{equation}

\paragraph{Visit selected positions}
\begin{equation}
\mbox{\begin{tabular*}{\dimexpr\linewidth-5em\relax}{@{\extracolsep{\fill}}lr@{}}%
$\displaystyle \xpass{j} + \xppass{j} + \sum_{i' \in \Ical_j : i' \geq i} \qxji{x}{ji'} + \sum_{i' \in \Ical_j : i' \leq i} \pxji{x}{ji'} \geq x_{ji}$ &
$\displaystyle j \in \Jcal,\; i \in \Ical_j$%
\end{tabular*}}
\end{equation}

\paragraph{Position selection binary}
\begin{equation}
\mbox{\begin{tabular*}{\dimexpr\linewidth-5em\relax}{@{\extracolsep{\fill}}lr@{}}%
$\displaystyle x_{ji} \in \{0,1\}$ &
$\displaystyle j \in \Jcal,\; i \in \Ical_j$%
\end{tabular*}}
\end{equation}

\paragraph{Aisle connectivity}
\begin{equation}
\mbox{\begin{tabular*}{\dimexpr\linewidth-5em\relax}{@{\extracolsep{\fill}}lr@{}}%
$\displaystyle \tilde{x}_j \geq x_{ji}$ &
$\displaystyle j \in \Jcal,\; i \in \Ical_j$%
\end{tabular*}}
\end{equation}
\begin{equation}
\mbox{\begin{tabular*}{\dimexpr\linewidth-5em\relax}{@{\extracolsep{\fill}}l@{}}%
$\displaystyle \tilde{x}_l = 1$%
\end{tabular*}}
\end{equation}
\begin{equation}
\mbox{\begin{tabular*}{\dimexpr\linewidth-5em\relax}{@{\extracolsep{\fill}}lr@{}}%
$\displaystyle \xx{x}_j + \xb{x}_j + \xbx{x}_j + \xxb{x}_j = \tilde{x}_{j+1}$ &
$\displaystyle j \in \Jcal \setminus \{m-1\} : j \geq l$%
\end{tabular*}}
\end{equation}
\begin{equation}
\mbox{\begin{tabular*}{\dimexpr\linewidth-5em\relax}{@{\extracolsep{\fill}}lr@{}}%
$\displaystyle \xx{x}_j + \xb{x}_j + \xbx{x}_j + \xxb{x}_j = \tilde{x}_j$ &
$\displaystyle j \in \Jcal : j < l$%
\end{tabular*}}
\end{equation}
\begin{equation}
\mbox{\begin{tabular*}{\dimexpr\linewidth-5em\relax}{@{\extracolsep{\fill}}lr@{}}%
$\displaystyle \tilde{x}_j \geq \tilde{x}_{j+1}$ &
$\displaystyle j \in \Jcal \setminus \{m-1\} : j \geq l$%
\end{tabular*}}
\end{equation}
\begin{equation}
\mbox{\begin{tabular*}{\dimexpr\linewidth-5em\relax}{@{\extracolsep{\fill}}lr@{}}%
$\displaystyle \tilde{x}_j \leq \tilde{x}_{j+1}$ &
$\displaystyle j \in \Jcal : j < l$%
\end{tabular*}}
\end{equation}
\begin{equation}
\mbox{\begin{tabular*}{\dimexpr\linewidth-5em\relax}{@{\extracolsep{\fill}}lr@{}}%
$\displaystyle \tilde{x}_j \in \{0,1\}$ &
$\displaystyle j \in \Jcal$%
\end{tabular*}}
\end{equation}

%% file: main.bbl
\begin{thebibliography}{27}
\expandafter\ifx\csname natexlab\endcsname\relax\def\natexlab#1{#1}\fi
\providecommand{\url}[1]{\texttt{#1}}
\providecommand{\href}[2]{#2}
\providecommand{\path}[1]{#1}
\providecommand{\DOIprefix}{doi:}
\providecommand{\ArXivprefix}{arXiv:}
\providecommand{\URLprefix}{URL: }
\providecommand{\Pubmedprefix}{pmid:}
\providecommand{\doi}[1]{\href{http://dx.doi.org/#1}{\path{#1}}}
\providecommand{\Pubmed}[1]{\href{pmid:#1}{\path{#1}}}
\providecommand{\bibinfo}[2]{#2}
\ifx\xfnm\relax \def\xfnm[#1]{\unskip,\space#1}\fi
\bibitem[{Boysen and De~Koster(2025)}]{boysen202550}
\bibinfo{author}{Boysen, N.}, \bibinfo{author}{De~Koster, R.},
  \bibinfo{year}{2025}.
\newblock \bibinfo{title}{50 years of warehousing research—an operations
  research perspective}.
\newblock \bibinfo{journal}{European Journal of Operational Research}
  \bibinfo{volume}{320}, \bibinfo{pages}{449--464}.
\bibitem[{Boysen et~al.(2019)Boysen, De~Koster and
  Weidinger}]{boysen2019warehousing}
\bibinfo{author}{Boysen, N.}, \bibinfo{author}{De~Koster, R.},
  \bibinfo{author}{Weidinger, F.}, \bibinfo{year}{2019}.
\newblock \bibinfo{title}{Warehousing in the e-commerce era: A survey}.
\newblock \bibinfo{journal}{European Journal of Operational Research}
  \bibinfo{volume}{277}, \bibinfo{pages}{396--411}.
\bibitem[{Daniels et~al.(1998)Daniels, Rummel and Schantz}]{daniels1998model}
\bibinfo{author}{Daniels, R.L.}, \bibinfo{author}{Rummel, J.L.},
  \bibinfo{author}{Schantz, R.}, \bibinfo{year}{1998}.
\newblock \bibinfo{title}{A model for warehouse order picking}.
\newblock \bibinfo{journal}{European Journal of Operational Research}
  \bibinfo{volume}{105}, \bibinfo{pages}{1--17}.
\bibitem[{De~Koster et~al.(2007)De~Koster, Le-Duc and
  Roodbergen}]{de2007design}
\bibinfo{author}{De~Koster, R.}, \bibinfo{author}{Le-Duc, T.},
  \bibinfo{author}{Roodbergen, K.J.}, \bibinfo{year}{2007}.
\newblock \bibinfo{title}{Design and control of warehouse order picking: A
  literature review}.
\newblock \bibinfo{journal}{European journal of operational research}
  \bibinfo{volume}{182}, \bibinfo{pages}{481--501}.
\bibitem[{Dunn et~al.(2025a)Dunn, Charkhgard, Eshragh and
  Stojanovski}]{dunn2025double}
\bibinfo{author}{Dunn, G.}, \bibinfo{author}{Charkhgard, H.},
  \bibinfo{author}{Eshragh, A.}, \bibinfo{author}{Stojanovski, E.},
  \bibinfo{year}{2025}a.
\newblock \bibinfo{title}{Double traversals in optimal picker routes for
  warehouses with multiple blocks}.
\newblock \bibinfo{journal}{Operations Research Letters} ,
  \bibinfo{pages}{107397}.
\bibitem[{Dunn et~al.(2025b)Dunn, Stojanovski, Lamichhane, Charkhgard and
  Eshragh}]{dunn2025deterministic}
\bibinfo{author}{Dunn, G.}, \bibinfo{author}{Stojanovski, E.},
  \bibinfo{author}{Lamichhane, B.}, \bibinfo{author}{Charkhgard, H.},
  \bibinfo{author}{Eshragh, A.}, \bibinfo{year}{2025}b.
\newblock \bibinfo{title}{Deterministic structure of vertical configurations in
  minimal picker tours for rectangular warehouses}.
\newblock \bibinfo{note}{Manuscript under review}.
\bibitem[{Goeke and Schneider(2021)}]{goeke2021modeling}
\bibinfo{author}{Goeke, D.}, \bibinfo{author}{Schneider, M.},
  \bibinfo{year}{2021}.
\newblock \bibinfo{title}{Modeling single-picker routing problems in classical
  and modern warehouses: Informs journal on computing meritorious paper
  awardee}.
\newblock \bibinfo{journal}{INFORMS Journal on Computing} \bibinfo{volume}{33},
  \bibinfo{pages}{436--451}.
\bibitem[{Gu et~al.(2007)Gu, Goetschalckx and McGinnis}]{gu2007research}
\bibinfo{author}{Gu, J.}, \bibinfo{author}{Goetschalckx, M.},
  \bibinfo{author}{McGinnis, L.F.}, \bibinfo{year}{2007}.
\newblock \bibinfo{title}{Research on warehouse operation: A comprehensive
  review}.
\newblock \bibinfo{journal}{European journal of operational research}
  \bibinfo{volume}{177}, \bibinfo{pages}{1--21}.
\bibitem[{{Gurobi Optimization, LLC}(2024)}]{gurobi2024}
\bibinfo{author}{{Gurobi Optimization, LLC}}, \bibinfo{year}{2024}.
\newblock \bibinfo{title}{Gurobi optimizer reference manual}.
\newblock \URLprefix \url{https://www.gurobi.com/documentation/}.
\bibitem[{Hall(1993)}]{hall1993distance}
\bibinfo{author}{Hall, R.W.}, \bibinfo{year}{1993}.
\newblock \bibinfo{title}{Distance approximations for routing manual pickers in
  a warehouse}.
\newblock \bibinfo{journal}{IIE transactions} \bibinfo{volume}{25},
  \bibinfo{pages}{76--87}.
\bibitem[{He{\ss}ler and Irnich(2022)}]{hessler2022note}
\bibinfo{author}{He{\ss}ler, K.}, \bibinfo{author}{Irnich, S.},
  \bibinfo{year}{2022}.
\newblock \bibinfo{title}{A note on the linearity of ratliff and rosenthal's
  algorithm for optimal picker routing}.
\newblock \bibinfo{journal}{Operations Research Letters} \bibinfo{volume}{50},
  \bibinfo{pages}{155--159}.
\bibitem[{He{\ss}ler and Irnich(2024)}]{hessler2024exact}
\bibinfo{author}{He{\ss}ler, K.}, \bibinfo{author}{Irnich, S.},
  \bibinfo{year}{2024}.
\newblock \bibinfo{title}{Exact solution of the single-picker routing problem
  with scattered storage}.
\newblock \bibinfo{journal}{INFORMS Journal on Computing} \bibinfo{volume}{36},
  \bibinfo{pages}{1417--1435}.
\bibitem[{Irnich and L{\"u}ke(2025)}]{irnich2025single}
\bibinfo{author}{Irnich, S.}, \bibinfo{author}{L{\"u}ke, L.},
  \bibinfo{year}{2025}.
\newblock \bibinfo{title}{The Single Picker Routing Problem with Scattered
  Storage in Parallel-Aisle Warehouse with Multiple Blocks}.
\newblock \bibinfo{type}{Technical Report}. Johannes Gutenberg-Universit{\"a}t
  Mainz.
\bibitem[{L{\"u}ke et~al.(2025)L{\"u}ke, Hessenius and Irnich}]{luke2025linear}
\bibinfo{author}{L{\"u}ke, L.}, \bibinfo{author}{Hessenius, A.},
  \bibinfo{author}{Irnich, S.}, \bibinfo{year}{2025}.
\newblock \bibinfo{title}{A linear-size model for the single picker routing
  problem with scattered storage}.
\newblock \bibinfo{journal}{European Journal of Operational Research} .
\bibitem[{L{\"u}ke et~al.(2024)L{\"u}ke, He{\ss}ler and
  Irnich}]{luke2024single}
\bibinfo{author}{L{\"u}ke, L.}, \bibinfo{author}{He{\ss}ler, K.},
  \bibinfo{author}{Irnich, S.}, \bibinfo{year}{2024}.
\newblock \bibinfo{title}{The single picker routing problem with scattered
  storage: modeling and evaluation of routing and storage policies}.
\newblock \bibinfo{journal}{OR Spectrum} \bibinfo{volume}{46},
  \bibinfo{pages}{909--951}.
\bibitem[{Pansart et~al.(2018)Pansart, Catusse and
  Cambazard}]{pansart2018exact}
\bibinfo{author}{Pansart, L.}, \bibinfo{author}{Catusse, N.},
  \bibinfo{author}{Cambazard, H.}, \bibinfo{year}{2018}.
\newblock \bibinfo{title}{Exact algorithms for the order picking problem}.
\newblock \bibinfo{journal}{Computers \& Operations Research}
  \bibinfo{volume}{100}, \bibinfo{pages}{117--127}.
\bibitem[{Petersen(1997)}]{petersen1997evaluation}
\bibinfo{author}{Petersen, C.G.}, \bibinfo{year}{1997}.
\newblock \bibinfo{title}{An evaluation of order picking routeing policies}.
\newblock \bibinfo{journal}{International Journal of Operations \& Production
  Management} .
\bibitem[{Prunet et~al.(2025)Prunet, Absi and Cattaruzza}]{prunet2025note}
\bibinfo{author}{Prunet, T.}, \bibinfo{author}{Absi, N.},
  \bibinfo{author}{Cattaruzza, D.}, \bibinfo{year}{2025}.
\newblock \bibinfo{title}{A note on the complexity of the picker routing
  problem in multi-block warehouses and related problems}.
\newblock \bibinfo{journal}{Annals of Operations Research}
  \bibinfo{volume}{347}, \bibinfo{pages}{1595--1605}.
\bibitem[{Ratliff and Rosenthal(1983)}]{ratliff1983order}
\bibinfo{author}{Ratliff, H.D.}, \bibinfo{author}{Rosenthal, A.S.},
  \bibinfo{year}{1983}.
\newblock \bibinfo{title}{Order-picking in a rectangular warehouse: a solvable
  case of the traveling salesman problem}.
\newblock \bibinfo{journal}{Operations research} \bibinfo{volume}{31},
  \bibinfo{pages}{507--521}.
\bibitem[{Roodbergen and De~Koster(2001)}]{roodbergen2001routing}
\bibinfo{author}{Roodbergen, K.J.}, \bibinfo{author}{De~Koster, R.},
  \bibinfo{year}{2001}.
\newblock \bibinfo{title}{Routing order pickers in a warehouse with a middle
  aisle}.
\newblock \bibinfo{journal}{European Journal of Operational Research}
  \bibinfo{volume}{133}, \bibinfo{pages}{32--43}.
\bibitem[{Rouwenhorst et~al.(2000)Rouwenhorst, Reuter, Stockrahm, van Houtum,
  Mantel and Zijm}]{rouwenhorst2000warehouse}
\bibinfo{author}{Rouwenhorst, B.}, \bibinfo{author}{Reuter, B.},
  \bibinfo{author}{Stockrahm, V.}, \bibinfo{author}{van Houtum, G.J.},
  \bibinfo{author}{Mantel, R.}, \bibinfo{author}{Zijm, W.H.},
  \bibinfo{year}{2000}.
\newblock \bibinfo{title}{Warehouse design and control: Framework and
  literature review}.
\newblock \bibinfo{journal}{European journal of operational research}
  \bibinfo{volume}{122}, \bibinfo{pages}{515--533}.
\bibitem[{Saylam et~al.(2024)Saylam, {\c{C}}elik and S{\"u}ral}]{saylam2024arc}
\bibinfo{author}{Saylam, S.}, \bibinfo{author}{{\c{C}}elik, M.},
  \bibinfo{author}{S{\"u}ral, H.}, \bibinfo{year}{2024}.
\newblock \bibinfo{title}{Arc routing based compact formulations for picker
  routing in single and two block parallel aisle warehouses}.
\newblock \bibinfo{journal}{European Journal of Operational Research}
  \bibinfo{volume}{313}, \bibinfo{pages}{225--240}.
\bibitem[{Tompkins et~al.(2010)Tompkins, White, Bozer and
  Tanchoco}]{tompkins2010facilities}
\bibinfo{author}{Tompkins, J.A.}, \bibinfo{author}{White, J.A.},
  \bibinfo{author}{Bozer, Y.A.}, \bibinfo{author}{Tanchoco, J.M.A.},
  \bibinfo{year}{2010}.
\newblock \bibinfo{title}{Facilities planning}.
\newblock \bibinfo{publisher}{John Wiley \& Sons}.
\bibitem[{Valle et~al.(2017)Valle, Beasley and Da~Cunha}]{valle2017optimally}
\bibinfo{author}{Valle, C.A.}, \bibinfo{author}{Beasley, J.E.},
  \bibinfo{author}{Da~Cunha, A.S.}, \bibinfo{year}{2017}.
\newblock \bibinfo{title}{Optimally solving the joint order batching and picker
  routing problem}.
\newblock \bibinfo{journal}{European journal of operational research}
  \bibinfo{volume}{262}, \bibinfo{pages}{817--834}.
\bibitem[{Weidinger(2018)}]{weidinger2018picker}
\bibinfo{author}{Weidinger, F.}, \bibinfo{year}{2018}.
\newblock \bibinfo{title}{Picker routing in rectangular mixed shelves
  warehouses}.
\newblock \bibinfo{journal}{Computers \& Operations Research}
  \bibinfo{volume}{95}, \bibinfo{pages}{139--150}.
\bibitem[{Weidinger et~al.(2019)Weidinger, Boysen and
  Schneider}]{weidinger2019picker}
\bibinfo{author}{Weidinger, F.}, \bibinfo{author}{Boysen, N.},
  \bibinfo{author}{Schneider, M.}, \bibinfo{year}{2019}.
\newblock \bibinfo{title}{Picker routing in the mixed-shelves warehouses of
  e-commerce retailers}.
\newblock \bibinfo{journal}{European Journal of Operational Research}
  \bibinfo{volume}{274}, \bibinfo{pages}{501--515}.
\bibitem[{Wildt et~al.(2025)Wildt, Weidinger and Boysen}]{wildt2025picker}
\bibinfo{author}{Wildt, C.}, \bibinfo{author}{Weidinger, F.},
  \bibinfo{author}{Boysen, N.}, \bibinfo{year}{2025}.
\newblock \bibinfo{title}{Picker routing in scattered storage warehouses: an
  evaluation of solution methods based on tsp transformations}.
\newblock \bibinfo{journal}{OR Spectrum} \bibinfo{volume}{47},
  \bibinfo{pages}{35--66}.

\end{thebibliography}
